\documentclass[11pt]{amsart}

\usepackage{latexsym}
\usepackage{amssymb}
\usepackage{color}
\usepackage[dvips,frame,cmtip,arrow,matrix,line,graph]{xy}

\input xy
\xyoption{all}

\usepackage{amsmath}
\usepackage{amsthm}
\usepackage{amsfonts}
\usepackage{amscd}

\usepackage{epsfig}

\usepackage{caption2}
\usepackage{mathrsfs}

\newtheorem{thm}{Theorem}[section]
\newtheorem{lem}[thm]{Lemma}
\newtheorem{prop}[thm]{Proposition}

\newtheorem{cor}[thm]{Corollary}
\newtheorem{Def}[thm]{Definition}

\newtheorem{Th}{Theorem}
\newtheorem{Co}[Th]{Corollary}

\newcommand{\vs}{\vspace{4mm}}

\newcommand{\A}{\mathcal{A}}

\newcommand{\BB}{\mathscr{B}}
\newcommand{\C}{\mathcal{C}}
\newcommand{\CC}{\mathbb{C}}
\newcommand{\CCC}{\mathscr{C}}
\newcommand{\D}{\mathcal{D}}

\newcommand{\F}{\mathcal{F}}
\newcommand{\G}{\mathcal{G}}
\newcommand{\GG}{\mathbb{G}}

\newcommand{\M}{\mathcal{M}}

\newcommand{\QQ}{\mathbb{Q}}
\newcommand{\RR}{\mathbb{R}}

\newcommand{\Z}{\mathbb{Z}}

\newcommand{\Dif}{\textrm{Diff}}

\newcommand{\al}{\alpha}

\newcommand{\Ga}{\Gamma}

\newcommand{\De}{\Delta}
\newcommand{\Da}{\overrightarrow{\Delta}}

\newcommand{\Om}{\Omega}
\newcommand{\s}{\sigma}
\newcommand{\Si}{\Sigma}
\newcommand{\pa}{\overrightarrow{p}\!}

\newcommand{\rar}{\longrightarrow}
\newcommand{\inc}{\hookrightarrow}
\newcommand{\sta}{\stackrel}
\newcommand{\minus}{\backslash}

\newcommand{\x}{\times}

\newcommand{\lgl}{\langle}
\newcommand{\rgl}{\rangle}

\newcommand{\del}{\partial}

\newcommand{\emp}{\varnothing}

\newcommand{\B}{\operatorname{B}\!}

\newcommand{\colim}{\operatorname{colim}}

\newcommand{\link}{\operatorname{Link}}

\newcommand{\RP}{\RR\textrm{P}^2}
\newcommand{\RPi}{\RR\textrm{P}^\infty}
\newcommand{\BO}{\operatorname{BO}}

\newcounter{samcounter}

\begin{document}

\bibliographystyle{plain}

\title[Stability for mapping class groups of non-orientable surfaces]{Homological stability for the mapping class groups of non-orientable surfaces}

\author{Nathalie Wahl}
 \address{Department of Mathematics\\
               University of Chicago\\
               5734 S. University Ave.\\
               Chicago IL 60637\\
               USA}
      \email{wahl@math.uchicago.edu}
      \thanks{Supported by a BAEF fellowship and NSF research grant DMS-05044932}

\date{\today}

\begin{abstract}
We prove that the homology of the mapping class groups of non-orientable surfaces stabilizes with the genus of the surface. 
Combining our result with recent work of Madsen and Weiss, we obtain that the classifying space of the stable mapping class group of non-orientable surfaces,
up to homology isomorphism, is the infinite loop space of a Thom spectrum built from the canonical bundle over the Grassmannians of 2-planes in $\RR^{n+2}$.
In particular, we show that the stable rational cohomology is a polynomial algebra on generators in degrees $4i$---this is the non-oriented analogue of the Mumford conjecture. 

\end{abstract}

\maketitle

\section{Introduction}
The mapping class group of a surface $S$ is the group of components of the space of diffeomorphisms of $S$. 
Since the components of the diffeomorphism group are contractible, except in a few cases \cite{EarEel,EarSch},
the classifying space of the mapping class group classifies surface bundles, and has therefore attracted great interest. 
The mapping class groups of orientable surfaces have been extensively studied and in particular, it was shown in the 1980's, in a 
fundamental paper by Harer \cite{Har85}, that its homology stabilizes 
with the genus of the surface. Moreover, the stable homology was computed recently \cite{Gal04}, 
using the breakthrough methods of Madsen and Weiss \cite{MadWei02}. 
A lot less is known about its non-orientable analogue despite its frequent appearance.

In this paper we prove that the homology of the mapping class groups of non-orientable surfaces stabilizes with the genus of the surface
 and we show that the stable homology is independent of the number of boundary components. 
Combining our result with the work of Madsen and Weiss, we give a computable homotopical model for the classifying space of the stable 
non-orientable mapping class group, up to homology isomorphism, built from Grassmannians of 2-planes in $\RR^{n+2}$. 
In particular, we obtain the non-oriented analogue of the Mumford conjecture, showing that the stable rational cohomology is a polynomial 
algebra $\QQ[\zeta_1,\zeta_2,\dots]$ with $\zeta_i$ in degree $4i$.

\vs

We now explain this in more detail. Let $S_{n,r}$  
denote a non-orientable surface of genus $n$ with $r$ boundary components, that is $S_{n,r}$ is a connected sum of $n$ copies of $\RP$ with 
$r$ discs removed. 
The {\em mapping class group} of $S_{n,r}$ is the group $\M_{n,r}=\pi_0{\rm Diff}(S_{n,r}\ {\rm rel}\ \del)$, 
the group of components of the diffeomorphisms of the surface which fix its boundary pointwise. 

When $r\ge 1$, there is a stabilization map $\al:\M_{n,r}\to \M_{n+1,r}$ 
obtained by gluing a M\"{o}bius band with a disc removed (or an $\RP$ with two discs removed)
to the surface and extending the diffeomorphisms by the identity on the added part.  
Similarly, there is a map  $\beta:\M_{n,r}\to \M_{n,r+1}$ obtained  by gluing a pair of pants to the surface. 
Gluing a disc moreover defines a map $\delta:\M_{n,r}\to \M_{n,r-1}$, which gives 
 a left inverse to $\beta$ when $r\ge 2$.

\begin{Th}[Stability Theorem]\label{main}
For any $r\ge 1$, 

{\rm (1)} $\al_i:H_i(\M_{n,r};\Z)\to H_i(\M_{n+1,r};\Z)$
is surjective when $n\ge 4i$ and an isomorphism when $n\ge 4i+3$.  

{\rm (2)} $\beta_i:H_i(\M_{n,r};\Z)\to H_i(\M_{n,r+1};\Z)$
is an isomorphism when $n\ge 4i+3$.  

{\rm (3)} $\delta_i:H_i(\M_{n,1};\Z)\to H_i(\M_{n,0};\Z)$ is is surjective when $n\ge 4i+1$ and an isomorphism when $n\ge 4i+5$. 
\end{Th}

The corresponding theorem  for mapping class 
groups of  orientable surfaces was proved in the 1980's by Harer \cite{Har85} and improved shortly after by Ivanov \cite{Iva90}.
They obtain respectively a slope 3 and slope 2 result,  
where the {\em slope} of a stability theorem with a linear bound is the slope of the bound.
Our stability theorem is with slope 4.
Note that an orientable genus corresponds to two non-orientable ones in the sense that
$(\#_g\, T^2)\ \#\ \RP\ \cong\ \#_{2g+1}\, \RP$
where $T^2$ denotes the torus. Hence  
our slope 4 corresponds to Ivanov's slope 2.

We prove our stability theorem using the action of the mapping class groups on complexes of arcs and circles in the surfaces. 
Two new phenomena occur in non-orientable surfaces:\\
- embedded circles can be 1-sided or 2-sided, and we define here a corresponding notion for arcs;\\ 
- arcs and circles may have orientable or non-orientable complement. \\
The main ingredient in proving Theorem~\ref{main} is the high connectivity of appropriate complexes of circles and 1-sided arcs. 
We define 1-sided arcs in any surface, orientable or not, equipped with a set of {\em oriented points} in its boundary.  
When a surface is orientable, the complex of all 1-sided arcs corresponds precisely to the complex 
of arcs between two sets of points used by Harer in \cite{Har85}. 
It turns out that there is a gap in Harer's proof of the connectivity of this complex, which we correct here. 
(This is our Theorem~\ref{FSD}, correcting \cite[Thm.~1.6]{Har85}, which is used in \cite[Thm.~2.12]{Iva90} so that both Harer and Ivanov rely 
on it for their stability result.) 
Our proof of (1) and (2) in Theorem~\ref{main} is logically independent of \cite{Har85} and \cite{Iva90}; 
we replace the train tracks in Harer's paper, and the Morse theory in Ivanov's, by a simple surgery technique of Hatcher \cite{Hat91}. 
We rely however on \cite{Iva90} for the connectivity of the complex of all circles in a surface to prove the last part of Theorem~\ref{main}.

Theorem~\ref{main} also holds verbatim for mapping class groups of non-orientable surfaces with a fixed number $s$ of punctures, which can be permutable or not. This can either been seen by inserting punctures everywhere in the proof of Theorem~\ref{main}, or as a corollary of Theorem~\ref{main} using a spectral sequence argument (see \cite{Han07}). 
A different stability theorem, for fixed genus but increasing number of permutable punctures, is proved in \cite{HatWah07}.

\vs

Let $S\to E\sta{p}{\to} B$ be a bundle of non-oriented surfaces. The {\em vertical tangent bundle} of $E$ is the 2-plane bundle $T_pE\to E$ which associates to 
each point of $E$ the tangent plane to the surface it lies in. This plane bundle has a first Pontrjagin class $p_1\in H^4(E)$. The powers of 
this class transferred down to $B$ define a family of classes $\zeta_1,\zeta_2,...\in H^*(B)$  with $\zeta_i$ in degree $4i$. 
The classes $\zeta_i$ are the non-oriented analogues of the Mumford-Morita-Miller classes \cite{Mil86,Mor87,Mum83}: 
When the surface bundle is orientable, the vertical tangent bundle 
has an Euler class $e\in H^2(E)$ which in turn produces a family of classes $\kappa_1,\kappa_2,...\in H^*(B)$ with $\kappa_i$ in degree $2i$. 
By construction, $\kappa_{2i}=\zeta_i$ when both classes are defined. 

The Mumford conjecture, as proved recently by Madsen and Weiss \cite{MadWei02}, says that the stable rational cohomology of the mapping class groups of 
orientable surfaces is the polynomial algebra $\QQ[\kappa_1,\kappa_2,\dots ]$. Madsen and Weiss actually prove a stronger result, 
using the reinterpretation of the classes \cite{MadTil01}, which  we now describe. 
Embedding $E$ in $B\times \RR^{n+2}$ for some large $n$ and collapsing everything outside a tubular neighborhood gives a map, the so-called Thom-Pontrjagin map,  
$S^{n+2}\wedge B\to Th(N_pE)$ to the Thom space of the vertical normal bundle of $E$. Classifying this bundle gives a further map into 
$Th(U_n^\perp)$, the Thom space of the orthogonal bundle to the tautological bundle over the Grassmannian $Gr(2,n+2)$ of $2$-planes in $\RR^{n+2}$ 
(or the Grassmannian of oriented planes if the bundle was oriented). 
The spaces $Th(U_n^\perp)$ form a spectrum  
 which we denote $\GG_{-2}$, and $\GG_{-2}^+=\CC P^\infty_{-1}$ in the oriented case. 
Applying the construction to the universal surface bundle gives a map 
$$\psi: \B\Dif(S)\rar \Om^\infty\GG_{-2}=\colim_{n\to \infty} \Om^{n+2}Th(U_{n}^\perp).$$ 
The rational cohomology of $\Om^\infty\GG_{-2}$ and $\Om^\infty\GG_{-2}^+$ are easily seen to be polynomial algebras, 
in the zeta and kappa classes respectively.

Except in a few low genus cases, the components of the diffeomorphism group $\Dif(S)$ of 
a surface $S$  (orientable or not) are contractible (\cite{EarEel,EarSch}, see also \cite{Gra73}), 
so that   $\B\Dif(S_{n,r})\simeq \B\M_{n,r}$ when $n+r$ is large enough, and similarly for orientable surfaces.  
The main result of \cite{MadWei02} is that the analogue of the above map $\psi$ in the oriented case 
induces an isomorphism between the homology of the stable mapping class group of orientable 
surfaces and the homology of $\Om^\infty\GG_{-2}^+$. 
This result relies on Harer's stability theorem. 
Forgetting orientations in \cite{MadWei02} and replacing Harer's result by our Theorem~\ref{main} yields our second main theorem. 
(We actually use the new simpler approach to the Madsen-Weiss theorem via cobordism categories \cite{GMTW05}.)

Let $\M_\infty$ denote the stable non-orientable mapping class group, that is $\M_\infty=\colim_{n\to \infty}\M_{n,1}$. 
Its homology is the stable homology of the non-orientable mapping class groups 
in the sense that $H_i(\M_{n,r})\cong H_i(\M_\infty)$ when $n\ge 4i+3$ if $r\ge 1$, or $n\ge 4i+5$ when $r=0$ by Theorem~\ref{main}.

\begin{Th}[Stable Homology]\label{next} The map $\psi$  induces an isomorphism\\ 
$H_*(\M_{\infty};\Z)\sta{\cong}{\rar} H_*(\Om_0^\infty \mathbb{G}_{-2};\Z)$. 
\end{Th}

Here $\Om^\infty_0\GG_{-2}$ denotes the 0th component of $\Om^\infty\GG_{-2}$. 
Let $Q$ denote the functor 
$\Om^\infty\Si^\infty=\colim_{n\to \infty} \Om^n\Si^n$ and $Q_0$ its 0th component. 
Looking away from $2$ or rationally, the right hand side in Theorem~\ref{next} simplifies further:

\begin{Co}  $H_*(\M_{\infty};\Z[\frac{1}{2}])\cong H_*(Q_0(\BO(2)_+);\Z[\frac{1}{2}])$
\end{Co}

\begin{Co}
$H^*(\M_\infty;\QQ)\cong \QQ[\zeta_1,\zeta_2,\dots]$ with $|\zeta_i|=4i$.
\end{Co}

The homology of $Q_0(\BO(2)_+)$ is known in the sense that 
$Q_0(\BO(2)_+)\simeq Q_0(S^0)\times Q(\BO(2))$ and  the homology 
$H_*(Q(X);\mathbb{F}_p)$ is an explicitly described algebra over $H_*(X;\mathbb{F}_p)$ for any prime $p$ \cite{DyeLas62,CohLadMay}. 
There is a fibration sequence of infinite loop spaces  
$$\Om^\infty \GG_{-2} \rar Q(\BO(2)_+) \rar \Om^\infty\RPi_{-1}$$
where $\RPi_{-1}=\GG_{-1}$ is the spectrum with $n$th space the Thom space of the orthogonal bundle to the canonical bundle
over $Gr(1,n)$. The homology of $\Om^\infty\RPi_{-1}$ is 2-torsion. 
The homology of $\Om^\infty_0\GG_{-2}$ can then be calculated using this fibration sequence as well 
as an analogous fibration sequence for $\Om^\infty\RPi_{-1}$ (see Section~\ref{homo}). The corresponding calculations in the orientable and spin case were 
carried out by Galatius \cite{Galspin,Gal04}.

Earlier calculations of the homology of $\M_{n,r}$ are sparse. 
Low dimensional homology calculations with $\Z/2$-coefficients where obtained by Zaw in \cite{Zaw04} for small genus and number of boundary components.
Korkmaz \cite{Kor98} and Stukow \cite{Stu07} calculated the first homology group and showed that  $H_1(M_{n,k})=\Z/2$ when $n\ge 7$.
This agrees with our result.

Let $\CCC_2$ denote the (1+1)-cobordism category with objects 
embedded circles in $\RR^\infty$ and morphisms the space of embedded (non-oriented) cobordisms between 
such circles. The morphism space between two sets of circles has the homotopy type of the disjoint union of the classifying spaces of the 
diffeomorphism groups of the cobordisms between the circles. 
Using Theorem~\ref{main}, we show that $\Om\B\CCC_2$ is homology equivalent to $\Z\x \B\M_\infty$. 
This is a non-oriented version of Tillmann's theorem \cite{Til97}. 
It says that the 
stable {\em non-orientable} mapping class group governs the homotopy type of the loop space of 
the cobordism category of all {\em non-oriented} surfaces. This reflects the fact that the 
non-orientable surface of arbitrarily high genus is a stable object for non-oriented surfaces in the sense that any surface can be ``stabilized'' 
into a non-orientable one by taking a connected sum with a copy of $\RP$. 
Theorem~\ref{next} then follows from the main theorem of \cite{GMTW05}, 
which says that $\B\CCC_2$ is a deloop of $\Om^\infty\GG_{-2}$.

\vs

Sections 2--4 in the paper give the  proof of  part (1) and (2) in Theorem~\ref{main}, that is the stability for surfaces with boundary. 
In Section 2, we prove the high connectivity of the complexes of 1-sided arcs and of arcs between two sets of points. In Section 3, we 
deduce such high connectivity results for the subcomplexes of arcs with non-orientable connected complement. 
In Section 4 we give the spectral sequence arguments. 
We prove (2) in Theorem \ref{main} by showing that $\beta$ induces an isomorphism on homology on the stable groups (Theorem~\ref{inf-stab}).
This is enough by the first part of Theorem~\ref{main} and has the advantage that one only needs to know that a stable complex is contractible, 
a weaker statement than that of the connectivity of a sequence of complexes tending to infinity. (A similar idea was used in \cite{HatVogWah}.) 
Section 5 proves the third part of Theorem~\ref{main} (given as Theorem~\ref{closedstab}) using the complex of circles with non-orientable 
connected complement. 
Finally Section 6 describes the stable homology. We first relate  the non-oriented (1+1)-cobordism category to the 
stable non-orientable mapping class group. We then use \cite{GMTW05} to  prove Theorem~\ref{next}.

I would like to thank Allen Hatcher, Nikolai Ivanov, Peter May, and in particular S\o ren Galatius for many helpful conversations. 
 I am also very grateful to Ib Madsen for getting me interested in this question. In fact, the paper started out as a joint project with him.

\section{1-sided arcs}

In this section, we consider two families of arc complexes: 1-sided arcs, defined below, and arcs from one set of points in the boundary of the surface 
to another. We use a surgery method from \cite{Hat91} to prove that the complexes are highly connected. 

\vs

Let $S$ be a surface, orientable or not. By an {\em arc in $S$}, we always mean an embedded arc intersecting $\del S$ only at 
its endpoints and doing so transversally. 
Isotopies of arcs fix their endpoints. We consider arcs with boundary on a fixed set of points $\De$ in $\del S$. 
 An arc in $(S,\De)$ is {\em non-trivial} if it is not isotopic to an arc of $\del S$ disjoint from $\De$.

Let $\Da$ be an {\em oriented set of points} in $\del S$, i.e.~each point of $\Da$ comes with a choice of an orientation 
of the boundary component of $S$ it lies in. Choosing an orientation for each boundary component of $S$, we thus have two types of points, 
those with positive orientation, i.e. the same as the chosen one, and those with negative orientation. If $S$ is orientable, a choice of 
orientation for $S$ gives a choice of 
orientation of its boundary components and thus a decomposition $\Da=\Delta^+\cup \De^-$, where each subset represents an orientation class of points. 
We assume that $(S,\Da)$ is {\em non-orientable}, i.e. that  either 
$S$ is non-orientable, or it has no orientation compatible with the orientations of the points of $\Da$ on its boundary. When $S$ is orientable, 
this means that  $\De^+$ and $\De^-$ must both be non-empty.

\begin{Def}
A {\em 1-sided arc} in $(S,\Da)$ is an arc in $S$ with boundary in $\Da$ and 
whose normal bundle can be oriented in a way which is compatible with the orientation of its endpoints.
\end{Def}

\begin{figure}
\includegraphics[width=.8\textwidth]{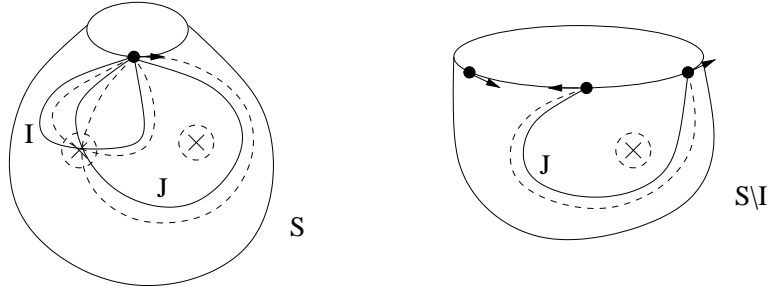}
\caption{1-sided arcs}\label{1-sided}
\end{figure}

Figure~\ref{1-sided} gives some examples of 1-sided arcs. (The crossed discs in the figure represent M\"{o}bius bands.) 
Note that an arc with boundary on a single point of $\Da$ is 1-sided if and only if it is a 1-sided curve. On the other hand, if $S$ is orientable, 
the 1-sided arcs are the arcs with a boundary point in $\Delta^+$ and the other one in $\Delta^-$. 

1-sided arcs as we defined them have the following property: A 1-sided arc stays 1-sided after cutting the surface along an arc whose interior is disjoint from it. As shown in Figure~\ref{1-sided}, 
this is only true when we keep track of the orientations of the endpoints of the arcs, which can be done by replacing the oriented points by  small oriented intervals and separating the arcs' endpoints on these intervals. In fact, a surface can become orientable after cutting 
along a 1-sided arc but the pair $(S,\Da)$ will stay non-orientable.

\begin{Def}
Let $S$ be a surface with non-empty boundary, $\De$, $\De_0$ and $\De_1$ be non-empty sets of points in $\del S$ with $\De_0\cap \De_1=\emp$, and 
let $\Da$ be a non-empty oriented set of points in $\del S$. 
We consider the following 3 simplicial complexes:\\
\begin{tabular}{lp{4in}}
$\A(S,\De)=$& the  vertices of $\A(S,\De)$ are isotopy classes of non-trivial embedded arcs with boundary points in $\De$. 
A $p$-simplex 
of $\A(S,\De)$ is a collection of $p+1$ such arcs $\langle I_0,\dots,I_p \rangle$ intersecting at most at their endpoints and not pairwise isotopic. 
\end{tabular}\\
\begin{tabular}{lp{4in}}
$\F(S,\Da)=$& the full subcomplex of $\A(S,\De)$ generated by the 1-sided arcs of $(S,\Da)$, where $\De$ is the underlying set of $\Da$.
\end{tabular}
\begin{tabular}{lp{3.65in}}
$\F(S,\De_0,\De_1)=$& the full subcomplex of $\A(S,\De_0\cup\De_1)$ of arcs with one boundary point in $\De_0$ and the other in $\De_1$.
\end{tabular}
\end{Def}

By the above remark, when $S$ is orientable, we have  
$$\F(S,\Da)\cong \F(S,\De^+,\De^-).$$ 
So $\F(S,\Da)$ and $\F(S,\De_0,\De_1)$ are two natural generalization of $\F(S,\De^+,\De^-)$ to non-orientable surfaces. 
We will see that the two complexes are closely related. 

\vs

Let $\De\subset \del S$ be the underlying set of $\Da$ or the union of two disjoint sets $\De_0$ and $\De_1$. 
We denote by $\del'S$ the components of $\del S$ intersecting $\De$ and if $r$ is the number of boundary components of $S$, $r'$ will denote 
the number of components of $\del' S$.  
The points of $\Da$ (resp.~$(\De_0,\De_1)$) subdivide the components of $\del'S$ into edges. 
We say that an edge in  $\del' S$ is {\em pure} if its endpoints have the same orientation (resp.~if its endpoints both belong to $\De_0$ or 
both to $\De_1$). An edge is {\em impure} if it is not pure. We say that $S$ has a {\em pure boundary component} if $\del'S$ has a component with 
only pure edges, i.e. where all the points are of the same type. 
Note that the number of impure edges is always even. We denote by $m$ the number of pure edges and by $l$ half the number of impure edges.

\begin{thm}\label{FSD}
{\rm (1)} If $\Da$ is a non-empty oriented set of points in $\del S$ such that the pair $(S,\Da)$ is non-orientable, then 
 $\F(S,\Da)$ is $(2h+r+r'+l+m-6)$-connected, 
where $h$ is the genus of $S$ if $S$ is non-orientable and twice its genus if $S$ is orientable,
and $r$ is the number of boundary components of $S$. 

{\rm (2)} If $\De_0$ and $\De_1$ are two disjoint non-empty sets of points in $\del S$, then $\F(S,\De_0,\De_1)$ is $(2h+r+r'+l+m-6)$-connected with $r$ and $h$ 
as above.  
\end{thm}

When $S$ is orientable, $\F(S;\De_0,\De_1)$  is the complex $BZ(\De,\De_0)$ of \cite{Har85} and  
the above result is a generalization of Theorem 1.6 in \cite{Har85}. There is however a gap in the proof: Harer's proof does not treat the case when 
$\De_0$ and $\De_1$ are two points in a single boundary component. 
We give here a complete and independent proof. 

To prove Theorem~\ref{FSD}, we need the following result: 

\begin{thm}\cite[Thm.~1.5]{Har85}\cite{Hat91}\label{ASD}
$\A(S,\De)$ is contractible unless $S$ is a disc or an annulus with $\De$ contained in one boundary component of $S$, in which case it is 
$(q+2r-7)$-connected, where $q$ is the order of $\De$ and $r=1,2$ is the number of boundary components of $S$. 
\end{thm}

This theorem was originally proved by Harer in \cite{Har85} for orientable surfaces using Thurston's theory of train tracks.
Hatcher reproved the result in \cite{Hat91} using a very simple surgery argument
which does not require the surface to be orientable. 
(Hatcher does not actually consider the case when $S$ is a disc or an annulus with $\De$ contained in a single boundary component, but this is 
easily obtained by induction from the cases $r=1$ with $q=4$ and $r=2$ with $q=2$ using an argument of the type used in Lemma~\ref{case3} below.)

We will need Theorem~\ref{ASD} in the last part of the argument for proving Theorem~\ref{FSD}, but we first 
 use the surgery argument of \cite{Hat91} to reduce the problem to a case which can be handled by a link argument.

\begin{figure}[ht]
\includegraphics[width=.9\textwidth]{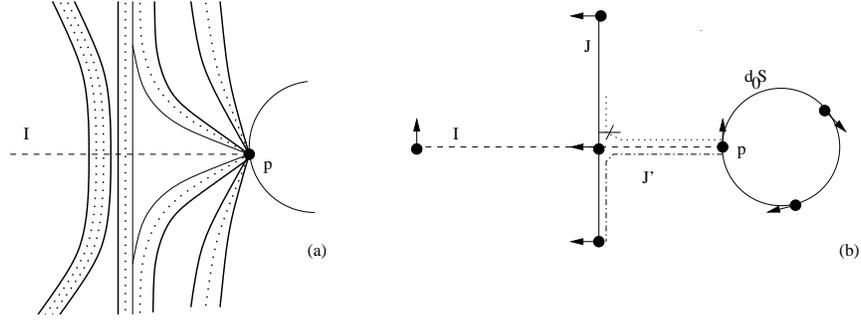}
\caption{Surgery for Lemma~\ref{case1}}\label{lem.case1}
\end{figure}

\begin{lem}\label{case1}
 If $S$ has at least one pure boundary component, then $\F(S,\Da)$ (resp. $\F(S,\De_0,\De_1)$) is contractible. 
\end{lem}

\begin{proof}
Let $\F$ denote $\F(S,\Da)$ or $\F(S,\De_0,\De_1)$ and
let $\del_0 S$ be a pure boundary of $S$. 
As $(S,\Da)$ is non-orientable (resp. $\De_0$ and $\De_1$ are both non-empty), 
there exists a non-trivial arc $I$ of $\F$ with at least one boundary at a point $p$ of $\del_0 S$. 
We define a deformation retraction of $\F$ onto the star of $\langle I\rangle$ by doing surgery along $I$ towards $p$, a process we describe now. 

Let $\s=\langle I_0,\dots,I_p \rangle$ be a simplex of $\F$ not already in the star of $I$. 
We can assume that the arcs $I_0,\dots,I_p$ intersect $I$ minimally, so that there 
is no disc whose boundary is the union of a subarc in $I$ and a subarc in some $I_j$. Any two such minimal positions for $\s$ are isotopic through minimal 
positions so that the following process will not depend on the choice of representing arcs for $\s$. 
Let $P=\sum_j t_jI_j$ be a point of $\s$ expressed in terms of barycentric coordinates. We think of $P$ as a weighted or thick family of 
arcs, where $I_j$ has thickness $t_j$. The arcs of $P$ cross $I$ with a total thickness $\theta=\sum_j a_jt_j$ where $a_j=|I\cap I_j|$. 
For $s\in [0,1]$, define $P_s$ to be the weighted family of arcs obtained from $P$ by cutting along $I$ from $p$ through the thickness up to 
$s\theta$ and sliding the newly created ends to $p$ as in Fig.~\ref{lem.case1}(a). For each arc we surger, we create two new arcs, one of which 
is in $\F$ and the other one, which we discard, is not. 
To see this in the first case, we consider an arc $J$ intersecting $I$ (see Fig.~\ref{lem.case1}(b)). As $J$ is 1-sided, its normal bundle can be oriented 
in a way which is compatible with its endpoints, so that it makes sense to consider the orientation of the normal bundle of $J$ at the intersection 
with $I$. Only one of the two arcs obtained by surgering $J$ along $I$ is 1-sided, namely $J'$ in the picture.
 Non-triviality of the resulting arc is insured by the fact that there is no trivial 1-sided arc with 
a boundary point on $\del_0 S$. (Either its other boundary is on a different boundary, or it is actually 1-sided as an arc in $S$ and hence non-trivial.)

For the second case, we just keep the only arc after surgery which still has a boundary in each of $\De_0$ and $\De_1$. 
Non-triviality follows from the fact that such an arc goes between different boundary components as one of its endpoints is in $\del_0 S$ which is pure. 

For any $s\in [0,1]$, $P_s$ is a point in a simplex of $\F$ and $P_s$ defines a path from $P$ to a point in a simplex in the star of $I$.
The key point here is that along the surgery process, when we pass from a simplex $\langle I_0,\dots,I_p \rangle$ to a simplex $\langle I'_0,\dots,I_p \rangle$, 
with $I_0'$ the non-trivial arc obtained by surgery along $I_0$, these two simplices are faces of a larger simplex of $\F$, namely 
$\langle I_0',I_0,\dots,I_p \rangle$. The path lies in this larger simplex. 
 The deformation respects the simplicial structure of $\F$ and hence defines a deformation of $\F$ onto the star of $I$ which is contractible. 
\end{proof}

Note that this surgery argument can be applied to many situations, as long as we make sure that surgery along any arc produces at least one non-trivial arc. 
This is what we will do in the next two lemmas. The last lemma uses a variant of the surgery argument. 

\begin{figure}[ht]
\includegraphics[width=.9\textwidth]{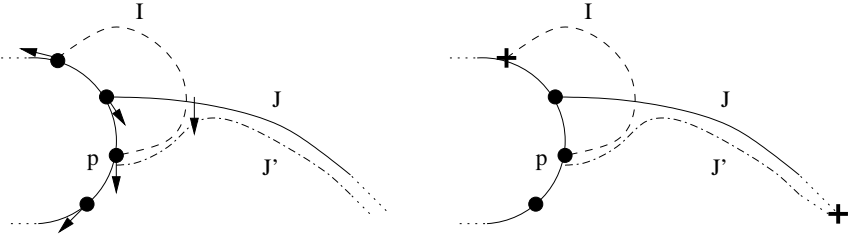}
\caption{Surgery for Lemma~\ref{case2}}\label{lem.case2}
\end{figure}

\begin{lem}\label{case2}
If $S$  has at least one pure  edge between a pure and an impure one in a  boundary component of $S$, 
then $\F(S,\Da)$ (resp.  $\F(S,\De_0,\De_1)$) is contractible. 
\end{lem}

\begin{proof}
As before, let $\F$ denote either $\F(S,\Da)$ or $\F(S,\De_0,\De_1)$. 
Consider the arc $I$ of $\F$ (drawn for both cases in Fig.~\ref{lem.case2}) 
which cuts a triangle with the impure and the first pure edge as the other sides of the triangle ---this arc is always non-trivial. We do surgery along $I$ towards $p$, i.e. towards 
the second pure edge, as in the previous lemma. Any arc $J$ which intersects $I$ minimally but non-trivially is as in Fig.~\ref{lem.case2}.
The surgery replaces $J$ by the arc $J'$ in the figure, which cannot be trivial because of the second pure edge. We thus get a 
deformation of $\F$ onto the star of $I$ as in Lemma~\ref{case1}.
\end{proof}

The following lemma is not strictly necessary for the proof of the theorem, but it will simplify the proof a little bit. 

\begin{figure}[ht]
\includegraphics[width=.8\textwidth]{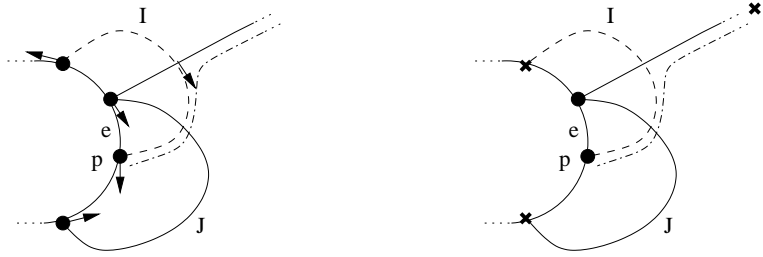}
\caption{Surgery for Lemma~\ref{case3}}\label{lem.case3}
\end{figure}

\begin{lem}\label{case3}
If the complex is non-empty, adding a pure edge between two impure edges increases the connectivity of $\F(S,\Da)$ (resp. $\F(S,\De_0,\De_1)$) by one. 
\end{lem}

\begin{proof}
Let $\F$ be the original complex ($\F(S,\Da)$ or $\F(S,\De_0,\De_1)$) and $\F'$ be the complex with an extra pure edge $e$ between two impure edges. 
There are two arcs $I$ and $J$ of $\F'$ (as shown in Fig.~\ref{lem.case3})  
cutting triangles with $e$ as one of the sides ---these arcs are non-trivial, unless $S$ is a disc and $\Delta$ has only 3 points, in which case the complex is empty. Note that the link of $\langle J\rangle$ is isomorphic to $\F$. 
If $\F'_J$ denotes the subcomplex of $\F'$ of simplices not containing $J$ as a vertex, we have 
$$\F'=\F'_J\bigcup_{{\rm Link}\langle J\rangle} {\rm Star}\langle J\rangle.$$ 
Now $\F'_J$ deformation retracts onto the star of $I$ by a surgery argument along $I$ towards $p$ as $J$ is the only arc that would 
become trivial after such a surgery. Hence $\F'$ 
is a suspension of the link of $\langle J\rangle$, and hence a suspension of $\F$. 
\end{proof}

For the next lemma, we restrict ourselves to the only case we really need for simplicity. 

\begin{figure}[ht]
\includegraphics[width=.9\textwidth]{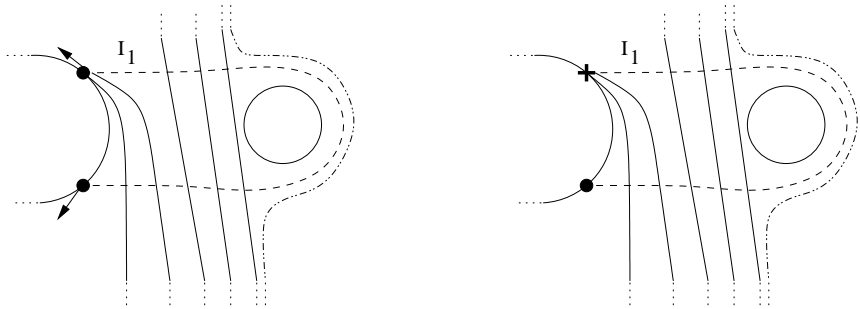}
\caption{Deformation for Lemma~\ref{case4}}\label{lem.case4}
\end{figure}

\begin{lem}\label{case4}
When $S$ has at least one impure edge and the complex is non-empty, adding a boundary component to $S$ disjoint from $\De$ increases the connectivity of 
$\F(S,\Da)$ (resp. $\F(S,\De_0,\De_1)$) by one. 
\end{lem}

\begin{proof}
Suppose that $S$ is obtained from $S'$ by capping off a boundary component $\del_0 S'$ disjoint from $\De$.
Let $\F$ be $\F(S,\Da)$ (resp. $\F(S,\De_0,\De_1)$) and   $\F'$ be $\F(S',\Da)$ (resp. $\F(S',\De_0,\De_1)$).  
We call a vertex $\langle I\rangle$ of $\F'$ {\em special} 
if it cuts an annulus with $\del_0 S'$ as one of its boundary and 
with the other boundary consisting of one impure edge of $\del S$ and a copy of $I$ itself. (See Fig.~\ref{lem.case4} 
where $I_1$ is a special vertex.) Special vertices exist when the complex is non-empty.
Let $\{I_j\}_{j\in \mathcal{J}}$ be the set of all the special vertices of $\F'$. Note that all the special vertices intersect, and hence 
$$\F'\ =\ \F'_{sp}\bigcup_{j\in \mathcal{J}} {\rm Star}\langle I_j\rangle$$
where $\F'_{sp}$ is the subcomplex of $\F'$ of simplices not containing the vertices $I_j$, $j\in \mathcal{J}$, and the union is taken along the links of the special vertices as $\F'_{sp}\cap{\rm Star}\langle I_j\rangle={\rm Link}\langle I_j\rangle$.  

We define a deformation retraction of $\F'_{sp}$ onto the star of $\langle I_1\rangle$ as follows.  
The intersection of a simplex with $I_1$ is always of the form shown in Fig.~\ref{lem.case4}. We deform the simplex as in the surgery argument 
but by passing the arcs one by one over $\del_0 S'$ instead of surgering them along some arc. 
The only arcs that would become trivial after such a deformation are the other special vertices as they would separate a trivial disc union 
an annulus with $\del_0S'$ as one of its boundaries, glued along an arc of $I_1$. Such an arc needs to have its two endpoints on the same boundary 
component of $S'$, with the endpoints forming an impure edge. 
We thus get a deformation retraction of $\F'_s$ onto the star of $I$.  
Noting now that the link of a special vertex is isomorphic to $\F$, we have that $\F'$ is a wedge of suspensions of copies of $\F$ and thus 
is one more connected. 
\end{proof}

\begin{proof}[Proof of Theorem~\ref{FSD}] (This part of the argument is missing in \cite{Har85}. It requires an intricate link argument 
which we learned from \cite{Iva87}.)
As before, let $\F$ denote either $\F(S,\Da)$ or $\F(S,\De_0,\De_1)$. 
The theorem is obviously true in the cases where we have already shown that $\F$ is contractible, i.e. when $S$ 
 has a pure boundary component (Lemma~\ref{case1}) or a component with an impure edge followed by at least two pure edges (Lemma~\ref{case2}). 
 We will prove the other cases by induction on the lexicographically ordered triple $(h,r,q)$, where the (doubled) genus $h\ge 0$, 
the number of boundaries $r\ge 1$ 
 and $q=2l+m\ge 1$ is the cardinality of $\Delta$ ($q\ge 2$ in the case $\F=\F(S,\De_0,\De_1)$). We can assume $l\ge 1$ and 
by Lemmas~\ref{case3}, \ref{case2} and \ref{case4}, it is enough to consider the case 
 when $m=0$  and $r=r'$ unless this case is empty. Emptyness occurs when $h=0$, $r=r'=1$ and $l=1,2$, in which case non-emptyness requires $r+l+m=4$. There are 5 such cases, which are easily checked, and the lemmas can be applied to these cases to increase $m$ and $(r-r')$.

So consider a surface $S$ and a set $\Da$ (resp. a pair of sets $(\De_0,\De_1)$) in $\del S$ with $r=r'$ and $m=0$ (and $l\ge 3$ if $h=0$ and $r=1$). 
We need to show that $\F$ is $(2h+2r+l-6)$-connected.
The induction starts with 
 $(h,r,q)=(0,1,6)$, where $q=2l$ as $m=0$. 
It is easy to check that $\F$ is non-empty in this case.

Fix now a surface $S$ with its corresponding tuple $(h,r,q)$ and let $k\le 2h+2r+l-6$ and 
 $f:S^k\to \F$ be any map. 
Forgetting the orientation of the points of $\Da$ (resp. taking $\De=\De_0\cup \De_1$) gives an inclusion $$\F\inc \A(S,\De).$$ 
When $\A(S,\De)$ is contractible, we can extend $f$ 
to a map of the disc in $\A(S,\De)$. To do this is general, by Theorem~\ref{ASD} we need $2h+2r+l-6\le 2r+q-7$ when $h=0$ and $r=1$ (as we assume $r'=r$), 
which is clear as $q=2l\ge 6$ in this case.

Choose an extension $\hat f$ of $f$:
$$\xymatrix{S^k \ar[r]^-f \ar@{^{(}->}[d] & \F \ar@{^{(}->}[d] \\
D^{k+1} \ar[r]^-{\hat f} & \A(S,\De)
}$$  
We may assume that $\hat f$ is simplicial with respect to some triangulation of $D^{k+1}$. 
We want to deform it so that its image lies in $\F$. We do this by induction.

Let $\s$ be a {\em bad} simplex of $D^{k+1}$, i.e. one whose image lies in $\A(S,\De)\minus \F$ so that all the arcs of $\s$ are pure.
 Cutting $S$ along the arcs of $\s$, we have in the first case 
$$(S,\Da)\minus \s=(X_1,\Da_1)\sqcup\dots\sqcup (X_c,\Da_c)\sqcup (Y_1,\Ga_1)\sqcup \dots\sqcup (Y_d,\Ga_d)$$ 
where each $(X_i,\Da_i)$ is non-orientable and each $Y_j$ is orientable in a way which is compatible with the points $\Ga_i$ inherited from $\Da$.   
(To obtain an orientation of the points of $\Da$ inherited in $S\minus \s$, replace each oriented point $p$ by a small oriented interval and pull apart  
the arcs with boundary at $p$, so that they have boundary on the oriented arc but do not intersect anymore. In $S\minus \s$, we now get small 
oriented intervals, which we can replace back by oriented points.)
In the second case, 
$$(S,\De_0,\De_1)\minus \s=(X_1,\De_0^1,\De^1_1)\sqcup\dots\sqcup (X_c,\De_0^c,\De_1^c)\sqcup (Y_1,\Ga_1)\sqcup \dots\sqcup (Y_d,\Ga_d)$$
where each $\De^i_\epsilon$ is a non-empty set in $\del S$ inherited from $\De_\epsilon$, and $\Ga_j$ is a subset of 
either $\De_0$ or $\De_1$, i.e. the $Y_j$'s have only points of one type. 

Let $Y_\s=i_1(Y_1)\cup \dots\cup i_d(Y_d)$, where $i_j:Y_j\to S$ is the canonical inclusion ---which is not necessarily injective on $\del Y_j$.  
Note that each component of $Y_\s$ has only points of one type in its boundary. 
We say that $\s$ is {\em regular} if no arc of $\s$ lies inside $Y_\s$. When $\s$ is regular, $Y_\s$ is the disjoint union of the $Y_j$'s with 
the points $p\in\Ga_j$ and $p'\in \Ga_{j'}$ identified if they come from the same point of $\De$. 

Let $\s$ be a regular bad $p$-simplex of $D^{k+1}$ maximal with respect to the ordered pair $(Y_\s,p)$, where 
$(Y',p')< (Y,p)$ if $Y'\varsubsetneq Y$ with $\del Y'\minus \del Y$ a union of non-trivial arcs in $Y$, or if $Y'=Y$ and $p'<p$. 
Then $\hat f$ restricts on the link of $\s$ to a map 
$${\rm Link}(\s)\cong S^{k-p} \to J_\s=\F(X_1)*\dots * \F(X_c)* \A(Y_1,\Ga_1)* \dots * \A(Y_d,\Ga_d)$$
where $\F(X_i)=\F(X_i,\Da_i)$ in the first case and   $\F(X_i)=\F(X_i,\De^i_0,\De_1^i)$  in the other case. 
Indeed, if a simplex $\tau$ in the link of $\s$ maps to pure arcs of $X_i$ for some $i$, then $Y_\s\subset Y_{\tau *\s}$ and either 
$Y_{\tau *\s}>Y_\s$, or $Y_{\tau *\s}=Y_s$ and $\s$ was not of maximal dimension. (On the other hand arcs of $\A(Y_j,\Ga_j)$ 
could not be added to $\s$ by regularity.)

In the case when one of the $Y_j$'s has non-zero genus or $b_j=b_j'>1$ boundaries, $\A(Y_j,\Ga_j)$ is contractible and hence so is $J_\s$. Suppose that all the 
$Y_j$'s have genus 0 and $b_j=1$, i.e. $Y_j$ is a disc. We want to calculate the connectivity of $J_\s$. 

We have $\chi(S\minus \s)=\chi(S)+p'+1$, where $p'+1\le p+1$ is the number of distinct arcs in the image of $\s$. 
If $X_i$ has genus $h_i$ and $r_i=r'_i$ boundaries, the above equation gives 

$\sum_{i=1}^c (2-h_i - r_i) + d= 2-h-r+p'+1$

$\Leftrightarrow \sum_{i=1}^c (h_i + r_i)=h+r-p'+2c+d-3.$\\
Moreover, we have 

$\sum_{i=1}^{c}m_i+\sum_{j=1}^d q_j=2p'+2$, where $m_i$ is the number of pure edges in $X_i$ and $q_j=|\Ga_j|$,  and $\sum_{i=1}^{c}l_i=l$,
where $l_i$ is half the number of impure edges in $X_i$.

By induction, $\F(X_i)$ is $(2h_i+2r_i+l_i+m_i-6)$-connected. Indeed for each $i$, we have 
$(h_i,r_i,q_i)<(h,r,q)$.  
On the other hand, $\A(Y_j,\Ga_j)$ is  $(q_j -5)$-connected. 
Thus  we get \\
${\rm Conn}(J_\s)\ge \sum_{i=1}^c (2h_i + 2r_i + l_i + m_i-4) + \sum_{j=1}^d (q_j -3)-2$

$=2h+2r-2p'+4c+2d-6+l+2p'+2-4c-3d-2$

$=2h+2r+l-d-6\ge 2h+2r+l-p-6\ge k-p$\\
because $3d\le p+1$, and hence $d\le p$, 
as the edges of $Y_j$ are arcs of $\s$, minimum three edges are needed for a non-trivial $Y_j$ and the same edge cannot be used for several 
$Y_j$'s by regularity of $s$.

So we can modify $\hat f$ on the interior of the star of $\s$ using a map 
$$\hat f * F\ :\ \del \s * D^{k-p+1}\cong {\rm Star}(\s) \ \rar\  \A(S,\De)$$ where $F:D^{k-p+1}\to J_\s$ is a map 
restricting to $\hat f$ on the boundary of $D^{k-p+1}$ which is identified with the link of $\s$. We are left to show that we have improved the 
situation this way. The new simplices are of the form $\tau=\al * \beta$ with $\al$ a proper face of $\s$ and $\beta$ mapping to $J_\s$. Suppose $\tau$ is 
a regular bad simplex. Then each arc of $\beta$ is pure and hence in $\A(Y_j,\Ga_j)$ for some $j$. Thus $Y_\tau\subset Y_\s$. If they are 
equal, we must have $\tau=\alpha$ is a face of $\s$ (by regularity of $\tau$ and $\s$). 
So $(Y_\tau,\dim(\tau))<(Y_\s,p)$ and the result follows by induction. 
\end{proof}

\section{Better 1-sided arcs}

To prove homological stability in the orientable case,  
one needs to pass to the subcomplexes of non-separating arc systems, 
that is arc systems  $\s=\langle I_0,\dots,I_p\rangle$ such that their complement $S\minus \s$ is connected. 
This serves the double purpose of reducing the number of orbits for the action of the mapping 
class group on these complexes, while restricting to simplices with the property that the inclusion of their stabilizer into the group is a composition of stabilization maps. When the surface is not orientable, in addition to discarding separating arc systems among the systems of 1-sided arcs, one also needs to eliminate arc systems with 
orientable complement. In this section, 
we use the results of Section 2 to prove that the complex $\G(S,\Da)$ of 1-sided arcs with connected non-orientable complement is highly connected. 
We also show that the analogous complex $\G(S,\pa_0,\pa_1)$ of 1-sided arcs between two points on different boundary components of the surface stabilizes to a contractible complex.

\begin{Def}
Let $S$ be any surface (orientable or not), $\Da$ an oriented set of points in $\del S$ and $\De_0,\De_1$ two disjoint sets of points in $\del S$.\\
\begin{tabular}{lp{4in}}
$BX(S,\Da)=$ & the subcomplex of $\F(S,\Da)$ of simplices with connected complement. 
\end{tabular}\\
\begin{tabular}{cp{3.55in}}
$BX(S,\De_0,\De_1)=$& the subcomplex of $\F(S,\De_0,\De_1)$ of simplices with connected complement. 
\end{tabular}

For $S$ non-orientable, we consider moreover the following subcomplexes:\\
\begin{tabular}{lp{4in}}
$\G(S,\Da)=$ & the subcomplex of $BX(S,\Da)$ of simplices with non-orientable (connected) complement. 
\end{tabular}\\
\begin{tabular}{lp{3.75in}}
$\G(S,\De_0,\De_1)=$& the subcomplex of $BX(S,\De_0,\De_1)$ of simplices with non-orientable (connected) complement. 
\end{tabular}
\end{Def}

One could also allow simplices of $\G(S,\Da)$ to have genus 0 complement, when all their faces are in $\G(S,\Da)$. However, this does not improve the stability range. (In our argument, this is reflected by the fact that the inequality in Proposition~\ref{G-action}(d) would not be improved.)  

\medskip

The following result was proved in the orientable case in \cite{Har85}. 
 We include the argument for convenience, even though it is standard.

\begin{thm}\label{BX}
$(1)$ If $\Da$ is a non-empty oriented set of points in $\del S$ such that $(S,\Da)$ is non-orientable, then $BX(S,\Da)$ is $(h+r'-3)$-connected, where $h$ is the genus of $S$ or twice its genus if $S$ is orientable, and $r'$ 
is the number of components of $\del S$ with points of $\Da$.

$(2)$ If $\De_0,\De_1$ are two disjoint non-empty sets of points in $\del S$, then the complex
 $BX(S,\De_0,\De_1)$ is $(h+r'-3)$-connected, for $h$ and $r'$ as above.
\end{thm}

\begin{proof}
We prove both cases in the theorem simultaneously by induction on the triple $(h,r,q)$, where $r$ is the number of components of $\del S$ and $q=|\Da|\ge 1$ (resp.~$q=|\De_0\cup \De_1|\ge 2$). 
To start the induction, note that the theorem is true when $h=0$ and $r'\le 2$ for any $r\ge r'$ and any $q$, and 
more generally that the complex is non-empty whenever $r'\ge 2$ or $h\ge 1$.

Let $BX$ denote either $BX(S,\Da)$ or $BX(S,\De_0,\De_1)$, and $\F$ denote $\F(S,\Da)$ or $\F(S,\De_0,\De_1)$ accordingly, and suppose that $(h,r,q)\ge(0,3,2)$ (as $q\ge 2$ when $h=0$). 
Then $h+r'-3\le 2h+r+r'+l+m-6$. Indeed, $r\ge 1$ and $l+m\ge 1$. Moreover we assumed that either $r\ge 3$ or $h\ge 1$. 

Fix $k\le h+r'-3$ and consider a map $f:S^k\to BX$. This map can be extended to a map $\hat f:D^{k+1}\to \F$ by 
Theorem~\ref{FSD}, with $\hat f$ simplicial for some triangulation of $D^{k+1}$. 
We call a simplex $\s$ of $D^{k+1}$ {\em regular bad} if $\hat f(\s)=\lgl I_0,\dots,I_p\rgl$ and each $I_j$ separates $S\minus\{I_0,\dots,\hat{I_j},\dots,I_p\}$. 
Let $\s$ be a regular bad simplex of maximal dimension. Suppose that $S\minus\s=X_1\sqcup\dots\sqcup X_c$. 
By maximality of $\s$, $\hat f$ restricts to a map 
$${\rm Link}(\s)\rar J_\s=BX(X_1)*\dots * BX(X_c)$$
where $BX(X_i)=BX(X_i,\Da_i)$ in the first case, with $(X_i,\Da_i)$ non-orientable as the arcs of $\s$ are 1-sided, and 
$BX(X_i)=BX(X_i,\De_0^i,\De_1^i)$ in the second case, with each $\De^i_\epsilon$  non-empty as the arcs of $\s$ are impure. Each $X_i$  has $(h_i,r_i,q_i)<(h,r,q)$, so by induction 
$BX(X_i)$ is $(h_i+r'_i-3)$-connected. The Euler characteristic gives $\sum_i(2-h_i-r_i')=2-h-n+p+1$. Now  $J_\s$ is 
$(\sum_i(h_i+r'_i-1)-2)$-connected, that is $(h+r'-p+c-5)$-connected. As $c\ge 2$, we can extend the restriction of $\hat f$ to $\link({\s})\simeq S^{k-p}$ to a map $F:D^{k+1-p}\to J_\s$. 
We modify $\hat f$ on the interior of the star of $\s$ using $\hat f * F$  on $\del \s * D^{k+1-p}\cong {\rm Star}(\s)$ as in the proof of Theorem~\ref{FSD}. 
If a simplex $\al * \beta$ in $\del\s * D^{k+1-p}$ is regular bad, $\beta$ must be trivial since simplices of $D^{k+1-p}$ do not separate $S\minus\hat f(\al)$, so that $\al * \beta=\al$ is a face of $\s$. We have thus  reduced the number of regular bad simplices of maximal dimension and the result follows by induction. 
\end{proof}

\begin{thm}\label{GSD}
Let $S$ be a non-orientable surface of genus $n\ge 1$. 

(1) For $\Da$ a non-empty oriented set of points in $\del S$, 
$\G(S,\Da)$ is $(n+r'-4)$-connected, where $r'$ is the number of components of $\del S$ intersecting $\Da$.

(2) For $\De_0,\De_1$ two non-empty disjoint sets of points in $\del S$, $\G(S,\De_0,\De_1)$ is $(n+r'-4)$-connected, with $r'$ as above. 
\end{thm}

\begin{proof}
Fix $k\le n+r'-4$ and let $\G$ denote either $\G(S,\Da)$ or $\G(S,\De_0,\De_1)$ and $BX$ denote $BX(S,\Da)$ or $BX(S,\De_0,\De_1)$ accordingly.
 Any map $f:S^k\to \G$ can be extended to a map $\hat f:D^{k+1}\to BX$ by Theorem~\ref{BX}, with $\hat f$ simplicial for some triangulation of $D^{k+1}$. 

A simplex $\s$ of $D^{k+1}$ is called {\em regular bad} if $\hat f(\s)=\lgl I_0,\dots,I_p\rgl$ with $S\minus\{I_0,\dots,I_p\}$ orientable but $S\minus\{I_0,\dots,\hat I_j,\dots,I_p\}$ non-orientable for each $j$.
 Note that each simplex $\lgl I_0,\dots,I_p\rgl$ of $BX$ with orientable complement contains a unique regular bad subsimplex, namely the subsimplex consisting of the arcs $I_j$ such that $S\minus\{I_0,\dots,\hat I_j,\dots,I_p\}$ is non-orientable. 

Let $\s$ be a regular bad simplex of maximal dimension in $D^{k+1}$. As $\s$ is maximal, $\hat f:\link(\s)\to BX(S\minus\hat f(\s))$. (Maximality is needed here in case there are vertices in the link of $\s$ mapping to some arc already in the image of $\s$.) The surface $S\minus\hat f(\s)$ is orientable of genus $g$ with $b$ boundary components, with $2-2g-b=2-n-r+p+1$. In particular, $BX(S\minus\hat f(\s))$ is $(n+r'-p-4)$-connected by Theorem~\ref{BX}. Hence we can extend the restriction of $\hat f$ to $\link(\s)\simeq S^{k-p}$ to a map $F:D^{k-p+1}\to BX(S\minus\hat f)$. We modify $\hat f$ in the interior of the star of $\s$ as in the previous theorems, using $\hat f * F:\del\s * D^{k+1-p}\cong {\rm Star}(\s) \to BX$. We claim that this reduces the number of regular bad simplices of maximal dimension. Indeed, consider a simplex $\al * \beta$  of $\del\s * D^{k+1-p}$. Suppose first that there exists an arc $I_j$ in $\hat f(\s)\minus\hat f(\al)$. 
Since $F(\beta)$ does not separate $S\minus\hat f(\s)$, there exists an arc $J$ in $S\minus(\hat f(\s) * F(\beta))$ joining the two copies of the midpoint of $I_j$. The arc $J$ closes up to a 1-sided circle in $S\minus(\hat f(\al)*F(\beta))$ as $S\minus \hat f(\s)$ is orientable and $S\minus\{I_0,\dots,\hat I_j,\dots,I_p\}$ is non-orientable. Hence $\hat f(\al) * F(\beta)$ has non-orientable complement. On the other hand, if $\hat f(\al)=\hat f(\s)$, then $\al$ is the regular bad subsimplex of $\al * \beta$, of dimension strictly smaller than $\s$. So modifying $\hat f$ in the interior of the star of $\s$ reduced the number of regular bad simplices of maximal dimension. The result follows by induction.
\end{proof}

\begin{Def}
Let $S$ be a non-orientable surface and $\Da_0$ and $\Da_1$ be two disjoint oriented sets of points in $\del S$. \\
\begin{tabular}{cp{3.65in}}
$\G(S,\Da_0,\Da_1)=$&$\G(S,\De_0,\De_1)\cap \G(S,\Da)$, where $\De_i$ is the underlying set of $\Da_i$ and $\Da=\Da_0\cup \Da_1$.  
\end{tabular}
\end{Def}

When $S$ is orientable, this corresponds to having four sets of points $\De_0,\De_1,\De_2,\De_3$ in $\del S$ and 
considering the complex of arcs between $\De_0$ and $\De_1$ and between $\De_2$ and $\De_3$ (with connected complement). These complexes can be shown to be highly connected in some particular cases 
but it seems hard to show high connectivity in general. 
For our purpose, it is however 
enough to prove a weaker statement, namely that the complex with $\Da_0$ and $\Da_1$ being single points in different 
boundary components of $S$ becomes contractible when the genus tends to infinity.

\begin{figure}[ht]
\includegraphics{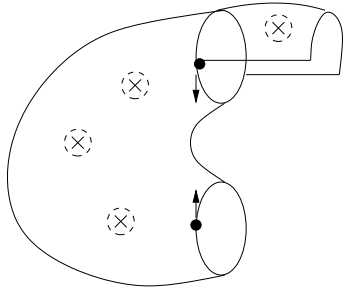}% 60%
\caption{The map $\alpha$}\label{alpha}
\end{figure}

Consider $\pa_0\in \del_0S$ and $\pa_1\in\del_1S$, two oriented points in a non-orientable surface $S$ with at least two boundary components $\del_0 S,\del_1 S$.
Let $$\al:S_{n,r}\to S_{n+1,r}$$ 
be the map that glues a M\"{o}bius band identifying part of its boundary to part of $\del_0 S$, away from $\pa_0$ (see Fig.~\ref{alpha}).
The map $\al$ induces an inclusion $\G(S_{n,r},\pa_0,\pa_1)\inc \G(S_{n+1,r},\pa_0,\pa_1)$. Let 
$$\G_{\infty,r}=\colim(\ \G(S_{n,r},\pa_0,\pa_1)\sta{\al}{\rar} \G(S_{n+1,r},\pa_0,\pa_1) \sta{\al}{\rar} \dots )$$

\begin{thm}\label{Ginf}
$\G_{\infty,r}$ is contractible for any $r\ge 2$. 
\end{thm}

\begin{proof}
Let $f:S^m\to \G_{\infty,r}$ be a map. As the sphere is compact, this map 
factors through $\G(S_{n,r},\pa_0,\pa_1)$ for some $n$ large. Now $\G(S_{n,r},\pa_0,\pa_1)\subset \G(S_{n,r},p_0,p_1)$ and we choose $n$ to be 
large enough so that the latter complex is at least $m$-connected. This is possible by Theorem~\ref{GSD}. Hence $f$ extends to a map 
$f':D^{m+1}\to \G(S_{n,r},p_0,p_1)$  which we can assume is simplicial for some triangulation of the disc. 
Our goal is to construct from $f'$ a map $f'':D^{m+1}\to\G(S_{n+2m+4,r},\pa_0,\pa_1)$ which extends $\al^{2m+4}\circ f$:
$$\xymatrix{S^m \ar[r]^-f \ar[d] & \G(S_{n,r},\pa_0,\pa_1) \ar[rr]^-{\alpha^{2m+4}} \ar[d] & & \G(S_{n+2m+4,r},\pa_0,\pa_1) \\
D^{m+1} \ar[r]^-{f'} \ar@/_3pc/@{-->}[urrr]_-{f''} & \G(S_{n,r},p_0,p_1)  & & }$$

\vs

For an arc $I\in\G(S_{n,r},p_0,p_1)$, there are $m+2$ canonical lifts to the complex $\G(S_{n+2m+4,r},\pa_0,\pa_1)$ defined as 
follows: we consider the $2m+4$ extra M\"{o}bius bands in $S_{n+2m+4,r}$ grouped two by two as in Fig.~\ref{lift}. For each $0\le j\le m+1$, define $I^j$ 
to be the arc which follows $I$ from $p_1$ to $p_0$ arriving on the left of $I$ at $p_0$ then goes through the $(2j+1)$st extra M\"{o}bius band, and through 
the $(2j+2)$nd if $I$ was already 1-sided, then back to $p_0$ as in the picture. Note that $\lgl I,I^0,\dots,I^m\rgl$ is a simplex 
of $\G(S_{n+2m+4,r},p_0,p_1)$ and of $\G(S_{n+2m+4,r},\pa_0,\pa_1)$ if $I$ was already 1-sided. 

For a simplex $\s=\lgl I_0,\dots,I_p\rgl$ of $\G(S_{n,r},p_0,p_1)$ with $p\le m+1$ and with  
 $I_0,\dots,I_p$  in the left to right order at $p_0$, we define the canonical lifting $\hat\s\in \G(S_{n+2m+4,r},\pa_0,\pa_1)$ by 
$$\hat\s=\lgl I^0_0,I_1^1,\dots,I^p_p\rgl.$$
(Fig.~\ref{lift} gives the case $p=1$ with $I_0$ a 2-sided arc and $I_1$ a 1-sided one.) 
The arcs have the property that $\lgl I_0^0,\dots,I_j^j,I_j,\dots,I_p\rgl$ is a $(p+1)$-simplex 
of $\G(S_{n+2m+4,r},p_0,p_1)$ for each $j$  
and of $\G(S_{n+2m+4,r},\pa_0,\pa_1)$ if the arcs $I_j,\dots,I_p$ are 1-sided. 
Moreover $\lgl I_0^{i_0},\dots,I_j^{i_j},I_j^{i_j'},I_{j+1}^{i'_{j+1}}\dots,I_p^{i_p'}\rgl$ is a simplex of $\G(S_{n+2m+4,r},\pa_0,\pa_1)$
whenever $i_0<\dots< i_j\le i_j'< i_{j+1}'< \dots< i_p'$.
\begin{figure}
\includegraphics{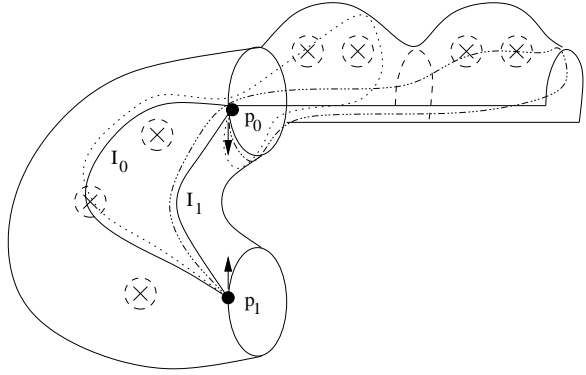} % 70% 
\caption{Lift of simplices}\label{lift}
\end{figure}

More generally, the arcs of a simplex $\s$ and their lifts form a poset defined by 
 $I<J$ for $I,J$ in $\s$ if $I$ is to the left of $J$ at $p_0$, $I^k<I$ for any $I$ in $\s$ and $I^k<J^l$ if $k<l$ and $I\le J$.  
A chain of inequalities in the poset corresponds to a simplex of $\G(S_{n+2m+4,r},p_0,p_1)$ and of 
$\G(S_{n+2m+4,r},\pa_0,\pa_1)$ if the arcs are 1-sided.

Let ${\rm Simp}(D^{m+1})$ denote the category of simplices of $D^{m+1}$. 
Consider the simplicial space (bisimplicial set) $$K={\rm hocolim}_{{\rm Simp}(D^{m+1})}F$$
where $F:{\rm Simp}(D^{m+1})\to {\rm SSets}$ is defined by $F(\s)=\Delta^p$ if $\s$ is a $p$-simplex of $D^{m+1}$.  
The realization of $K$ gives a new cell structure on $D^{m+1}\cong K$ with a cell $\De^p\times \De^q$ for each 
$p$-simplex $\tau_0$ of $D^{m+1}$ and $q$-chain of face inclusions $\tau_0\inc\dots\inc \tau_q$. 
We think of it as a cellular decomposition of the top simplices of $D^{m+1}$. 
(In Fig.~\ref{decomp} we show what a 2-simplex of $D^{m+1}$ is replaced by in $K$.) 

\begin{figure}[ht]
\input{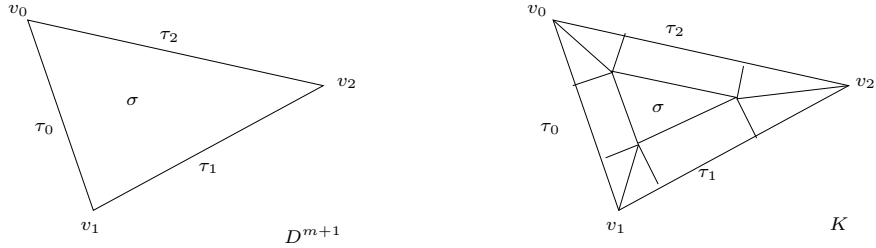} % 55%
\caption{$D^{m+1}$ and $K$}\label{decomp}
\end{figure}

We define $f'':K\to \G(S_{n+2m+4,r},\pa_0,\pa_1)$ as follows. First define $f''$ on the cells of the form $\De^p\times \{*\}$ associated to the 
simplices of $D^{m+1}$. If $\tau$ is a $p$-simplex of $D^{m+1}$, $f''(\tau)=f(\tau)$ if $\tau$ lies in 
the boundary of $D^{m+1}$, and $f''(\tau)=\hat \tau$ if $\tau$ lies in the interior of $D^{m+1}$. On the boundary, this is the original map $f$, 
and thus defines a simplicial map, but on the interior of $D^{m+1}$, if $\tau$ is a face of $\s$, $f''(\tau)$ is not necessarily a face of 
$f''(\s)$, as the lifts of the arcs of $\tau$ may be different in $\hat \tau$ and $\hat\s$. 
It is to make up for this difference that we have to replace $D^{m+1}$ by $K$.

For each chain of face inclusions $\bar \tau=\tau_0\inc\dots\inc \tau_q$ in $D^{m+1}$,   
with  $\tau_0=\lgl I_0<\dots<I_p\rgl$ in the interior of $D^{m+1}$,  
we define $f''(\bar\tau)$ to be the following $\De^p\times \De^q$ subcomplex of  $\G(S_{n+2m+4,r},\pa_0,\pa_1)$:
$$\begin{tabular}{ccccc}
$I_0^0$&$<$&$\dots$ &$<$&$I_p^p$ \\
$\wedge$&& $\dots$&&$\wedge$\\
$I_0^{i^1_0}$&$<$&$\dots$&$<$&$I_p^{i^1_p}$\\
$\wedge$&& $\dots$&&$\wedge$\\
$\vdots$&& $\vdots$&&$\vdots$\\
$\wedge$&& $\dots$&&$\wedge$\\
$I_0^{i^q_0}$&$<$&$\dots$&$<$&$I_p^{i^q_p}$
\end{tabular}$$
where $I_j^{i^k_j}$ is the lift of $I_j$ in $\hat\tau_k$. As $\tau_0$ is a face of $\tau_k$, we must have $i^k_j\ge j$. Similarly, $i^k_j\le i^l_j$ 
whenever $k<l$ as $\tau_k$ is a face of $\tau_l$.  Hence this forms a $\De^p\times \De^q$-subcomplex of $\G(S_{n+2m+4,r},\pa_0,\pa_1)$.
If $\tau_0, \dots, \tau_j$ are in the boundary of $D^{m+1}$, modify the definition of $f''$ above by using $f(\tau_0)$ on 
the $k$th row in the table for each $0\le k\le j$ and reversing some inequalities: 
$$\begin{tabular}{ccccc}
$I_0$&$<$&$\dots$ &$<$&$I_p$ \\
$\vee$&& $\dots$&&$\vee$\\
$I_0^{i^1_0}$&$<$&$\dots$&$<$&$I_p^{i^1_p}$
\end{tabular}$$
The map $f''$ respects the simplicial structure of $K$ and  is an extension of $f$ by construction. 
\end{proof}

\section{spectral sequence argument}

In this section, we prove the stability theorem for surfaces with boundaries.  We first show that the map $\al:\M_{n,r}\to \M_{n+1,r}$ induces an isomorphism 
in homology in a range, using the action of $\M_{n,r}$ on the complex $\G(S_{n,r},\pa)$ of {\em good} 1-sided arcs with boundary on a single point $p$.
 We then show that the stabilization map $\beta:\M_{n,r}\to\M_{n,r+1}$ induces an isomorphism in homology 
between the stable groups $\M_{\infty,r}\to \M_{\infty,r+1}$. For this, we use the action of $\M_{n,r+1}$ on the stable complex $\G_{\infty,r+1}$. 
Combining the two results gives the desired stability statement for $\beta$.

For the proof of (1) in Theorem \ref{main}, 
we use the relative spectral sequence argument.  
With a little extra work, one can use the non-relative argument, which in the orientable case gives a better stability bound (with slope 2 instead of slope 3). 
However, that second argument does not improve the bound in our case. 
The slope 4 comes from the geometry of the complexes we use, namely from the fact that cutting along an arc in a non-orientable surface can 
reduce the genus by 2 (that is by 1 orientable genus).

\subsection{Stabilization by projective planes}

Let $S_{n,r}$ be a non-orientable surface of genus $n$ with $r\ge 1$ boundary components, and consider $\G(S_{n,r},\Da)$ where 
$\Da=\overrightarrow{p}$ is a single point in a boundary component $\del_0 S$ of $S$. The mapping class group 
$\M_{n,r}:=\pi_0{\rm Diff}(S_{n,r}\ {\rm rel}\ \del)$ 
acts on $\G(S,\pa)$. 
The map $\al:S_{n,r}\to S_{n+1,r}$ shown in Fig.~\ref{alpha} induces an inclusion of simplicial complexes 
$$\al:X_n=\G(S_{n,r},\Da)\rar X_{n+1}=\G(S_{n+1,r},\Da)$$
which is $\M_{n,r}$-equivariant, where $\M_{n,r}$ acts on $X_{n+1}$ via the map $\M_{n,r}\to\M_{n+1,r}$ also induced by $\al$, extending 
the diffeomorphisms by the identity on the added M\"{o}bius band.

\begin{prop}\label{G-action}
 The action of $\M_{n,r}$ on $\G(S_{n,r},\pa)$ has the following properties:

{\rm (a)} $\M_{n,r}$ acts transitively on the vertices of $\G(S_{n,r},\pa)$.

{\rm (b)} An upper bound on the number of orbits of $q$-simplices is $\frac{(2q+2)!}{(q+1)!2^{q+1}}$, the number of possible orderings of the arcs at $\pa$.

{\rm (c)}  For a $q$-simplex $\s$ of $\G(S_{n,r},\pa)$, 
the possible diffeomorphism types of its complement $S\minus\s$ are $S_{n-k,r+k-q-1}$ where $q+1\le k\le 2q+1$. 
Moreover, the stabilizer of $\s$ is isomorphic to $\M(S\minus\s)$.

{\rm (d)} The map $\al:\G(S_{n,r},\pa)\to \G(S_{n+1,r},\pa)$ gives a 1--1 correspondence between the orbits of 
$q$-simplices of $\G(S_{n,r},\pa)$ and the orbits of $q$-simplices of $\G(S_{n+1,r},\pa)$ whenever $n\ge 2q+2$. 
\end{prop}

\begin{proof}
Let $\s$ be a $q$-simplex of $\G(S_{n,r},\pa)$. By assumption, $S\minus\s$ is connected and non-orientable. 
The ordering of the arcs of $\s$ at $\pa$, combined with their 1-sidedness in $S$, determines the number of boundary components of $S\minus \s$, and 
hence its genus by Euler characteristic. The ordering of the arcs also determines the {\em boundary pattern},
i.e. how the copies of the arcs of $\s$ lie in $\del(S\minus \s)$, 
so that if $\s$ and $\s'$ have the same ordering type at $\pa$, there is a diffeomorphism 
$S\minus\s\to S\minus\s'$ which takes the $j$th arc of $\s$ to the $\tau(j)$th arc of $\s'$ for some signed permutation $\tau\in \Si_{q+1}\int\Si_2$. 
This proves (b), and (a) follows from the fact that there is only one way to order one arc. 

Cutting along the arcs of $\s$ one by one,  the first arc is always 1-sided, so the 
total number of boundary components, which cannot decrease, can at most increase by $q$. 
Hence $S\minus \s$ has $r+j$ boundary components and genus $n-j-q-1$ for some $0\le j\le q$.  Moreover, a diffeomorphism 
of $S$ which fixes $\del S$ pointwise and fixes $\s$, fixes $\s$ pointwise, and is isotopic to a diffeomorphism fixing the arcs of $\s$ pointwise. This proves (c). 

If $n\ge 2q+2$, all possible orderings of the arcs of a $q$-simplex are realizable. Indeed, they can be obtained from a surface 
of genus $n-j-q-1$, with $j$ the same constant as above, by doing appropriate boundary identifications, these identifications forming the arcs 
of a simplex. This is possible precisely when $n-j-q-1\ge 1$ so that the complement of the arcs is a non-orientable surface.
As the orbits of $q$-simplices  are in 1--1 correspondence with the realizable orderings at $\pa$, there are in 1--1 correspondence for any $n\ge 2q+2$, and $\al$ induces such a correspondence.
\end{proof}

\begin{thm}\label{genusstab}
The map $\al: H_i(\M_{n,r};\Z)\rar H_i(\M_{n+1,r};\Z)$
is surjective when $n\ge 4i$ and is an isomorphism when $n\ge 4i+3$.  
\end{thm}

\begin{proof}
We use the relative spectral sequence argument of \cite{Vog79}, which we briefly recall here.

Let $X_n$ denote $\G(S_{n,r},\pa)$ as before. 
For short, we write $\M_n$ for $\M_{n,r}$. 
The map $\al$ induces a map of double complexes:
$$E_*\M_n \otimes_{\M_n} \tilde C_*(X_n) \rar E_*\M_{n+1} \otimes_{\M_{n+1}} \tilde C_*(X_{n+1})$$
where $E_*\M_n$ is a free resolution of $\M_n$ and $\tilde C_*(X_n)$ is the augmented chain complex of $X_n$. 
We take a level-wise cone of this map to obtain a double complex
$$C_{p,q}=(E_{q-1}\M_n \otimes_{\M_n} \tilde C_p(X_n))\oplus (E_q\M_{n+1} \otimes_{\M_{n+1}} \tilde C_p(X_{n+1}))$$
and consider the two associated spectral sequences.

The first spectral sequence, taking the homology in the $p$ direction first, has $E^1_{p,q}=0$ when $p\le n-3$ by 
Theorem~\ref{GSD} as $r'=1$.
The second spectral sequence thus converges to 0 when $p+q\le n-3$. Using Shapiro's lemma, 
when $n\ge 2p+2$ this second sequence has $E^1$-term 
$$E^1_{p,q}=\bigoplus_{\s\in\mathcal{O}_p}H_q(St_{n+1}(\s),St_{n}(\s))\, , \ \ \ E^1_{-1,q}=H_q(\M_{n+1},\M_n)$$ 
where $\mathcal{O}_p$, for $p\ge 0$, is the set of orbits of $p$-simplices in $X_n$ or $X_{n+1}$ (which are isomorphic under the assumption), 
$St_n(\s)$ denotes the stabilizer of $\s$ in $X_n$ and the inclusion $St_n(\s)\to St_{n+1}(\s)$ is induced by $\alpha$.

\vs

We prove the theorem by induction on the degree of the homology. The result is obviously true when $i=0$ for any $n,r$. Suppose it is 
true for all $i<q$. We start by showing surjectivity when $i=q$. So suppose $n\ge 4q$. 

First consider the boundary map $d^1:E^1_{0,q}=H_q(St_{n+1}(\s_0),St_n(\s_0))\to E^1_{-1,q}=H_q(\M_{n+1},\M_{n})$. 
As $n\ge 4q$ and $q\ge 1$, we have $q-1\le n-3$ and thus $E^{1}_{-1,q}$ must be killed before we reach $E^\infty$. 
The sources of differentials to $E^1_{-1,q}$ are the terms $E^1_{s,q-s}=\bigoplus_{\s\in\mathcal{O}_s}H_{q-s}(St_{n+1}(\s),St_n(\s))$, where 
$0\le s\le q$. (We need $n\ge 2s+2$ for the isomorphism between the sets of orbits, and this is satisfied as $s\le q$ and $q\ge 1$.) 
When $s\ge 1$, by induction we have $E_{s,q-s}^1=0$. Indeed, the stabilizer of an $s$-simplex $St_n(\s)\cong\M_{n-k,r+k-s-1}$ where 
$s+1\le k\le 2s+1$. By induction, the map $H_{q-s}(St_{n}(\s_s))\to H_{q-s}(St_{n+1}(\s_s))$ is surjective when $n-k\ge 4(q-s)$, 
which is satisfied if $n-2s-1\ge 4q-4s$, i.e. if $n\ge 4q-2s+1$, and we assumed $n\ge 4q$ and $s\ge 1$. The map 
$H_{q-s-1}(St_{n}(\s_s))\to H_{q-s-1}(St_{n+1}(\s_s))$ is injective when $n-k\ge 4(q-s-1)+3=4q-4s-1$, which is a weaker condition. 
So the relative groups are zero under our assumption and the map 
$d^1$ above must be surjective. Now a diagram chase in 
$$\xymatrix{H_q(St_{n+1}(\s_0))\ar[r]\ar[d] & H_q(\M_{n+1}) \ar[d]\\
H_q(St_{n+1}(\s_0),St_{n}(\s_0))\ar[d]_0\ar@{->>}[r]^-{d^1}  & H_q(\M_{n+1},\M_{n})\\
H_{q-1}(St_{n}(\s_0))\ar[d]^-\cong_-\al & \\
H_{q-1}(St_{n+1}(\s_0)) & 
}$$ 
shows that $ H_q(\M_{n+1},\M_{n})=0$. Here we use that $n-1\ge 4(q-1)+3$, i.e. $n\ge 4q$, so that the bottom map is an isomorphism, 
and that the composition $H_q(St_{n+1}(\s_0))\to H_q(\M_{n+1})\to  H_q(\M_{n+1},\M_{n})$ is zero because the following diagram commutes: 
$\xymatrix{ H_q(St_{n+1}(\s_0))\ar[r] \ar[d]^\cong &  H_q(\M_{n+1}) \ar[d]^{Id}\\ 
 H_q(\M_n) \ar[r]^\al & H_q(\M_{n+1}) 
}$
This is equivalent to the fact that $M_{n+1,r}$ acts transitively on the 0-simplices of $X_{n+1}$.

To show injectivity, we consider $d^1:E^1_{0,q+1}\to E^1_{-1,q+1}$. 
Assume now that  $n\ge 4q+3$. We have $q\le n-3$, so $E^1_{-1,q+1}$ must be killed. The sources of differentials to it are 
the terms $E^1_{s,q-s+1}=\bigoplus_{\s\in\mathcal{O}_s}H_{q-s+1}(St_{n+1}(\s),St_n(\s))$, where 
$0\le s\le q+1$ (as $n\ge 2q+4$). 
When $s\ge 1$, each of these is 0 by induction. Indeed, we need $n-k\ge 4(q-s+1)$ for all $s+1\le k\le 2s+1$, i.e. 
$n\ge 4q-2s+5$ which is satisfied by assumption as $s\ge 1$. (We also need $n-k\ge 4(q-s)+3$, which is a weaker condition.)
This shows that the above $d^1$ must be surjective. 
A diagram chase in 
 $$\xymatrix{ H_{q+1}(St_{n+1}(\s_0),St_{n}(\s_0)) \ar[r]\ar@{->>}[d]^{d^1} & H_{q}(St_{n}(\s_0)) \ar[d] & \\
 H_{q+1}(\M_{n+1},\M_{n}) \ar[r] & H_{q}(\M_{n}) \ar[r]^-\al & H_{q}(\M_{n+1}) 
}$$ 
shows that $\al$ must be injective. Here we use that there is an isomorphism $i$ making the following diagram commute:
$\xymatrix{H_{q}(St_{n}(\s_0)) \ar[r]^-\al \ar[d]^{id} & H_{q}(St_{n+1}(\s_0)) \ar[d]_i^\cong\\
H_{q}(St_{n}(\s_0)) \ar[r] & H_{q}(\M_n)}$
which shows that the composition $H_{q+1}(St_{n+1}(\s_0),St_{n}(\s_0)) \to H_{q}(St_{n}(\s_0))$ $ \to H_{q}(\M_n)$ is zero. 
(Thinking of $\M_n$ as $St_{n+1}(\s_0')$ and $St_n(\s_0)$ as $St_{n+1}(\s_0 * \s_0')$ for $\s_0'$ in $X_{n+1}$ defined by an arc in $S_{n+1,r}-S_{n,r}$,  this isomorphism $i$ can be obtained by conjugating by a diffeomorphism of $S_{n+1,r}$ which maps $\s_0$ to $\s_0'$ and fixes everything else outside a neighborhood of $\s_0*\s_0'$.)
\end{proof}

\subsection{Stabilization with respect to boundaries}

\begin{figure}[ht]
\input{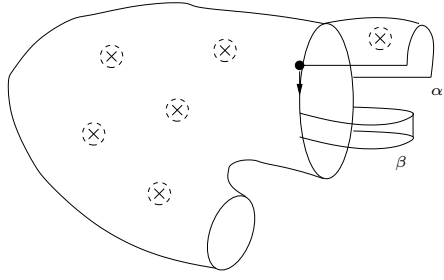} % 50%
\caption{$\beta$ commutes with $\al$}\label{mu}
\end{figure}

For $r\ge 1$, let $\beta:\M_{n,r}\to \M_{n,r+1}$ be the map induced by gluing a band on $\del_0S$ as in Fig.~\ref{mu}.  
As shown in the figure,  $\beta$ commutes with $\al:\M_{n,r}\to \M_{n+1,r}$, so that it induces a map 
on the stable groups $\beta:\M_{\infty,r}\to \M_{\infty,r+1}$, where  $\M_{\infty,r}=\colim(\M_{n,r} \sta{\al}{\to}\M_{n+1,r} \sta{\al}{\to}\dots)$.

Recall that $\G_{\infty,r}=\colim(\G(S_{n,r};\pa_0,\pa_1)\sta{\al}{\to} \G(S_{n+1,r};\pa_0,\pa_1)\sta{\al}{\to}\dots)$, 
where $\pa_0$ and $\pa_1$ are each single oriented points on different 
boundary components of $S$  and $r\ge 2$.
The action of $\M_{n,r}$ on $\G(S_{n,r};\pa_0,\pa_1)$ extends to an action of $\M_{\infty,r}$ on $\G_{\infty,r}$.

\begin{prop}\label{action2} The action of $\M_{\infty,r}$ on $\G_{\infty,r}$ has the following properties:

{\rm (a)} There is only one orbit of vertices.  

{\rm (b)} There are $(p+1)!$ orbits of $p$-simplices, corresponding to the possible changes of ordering from $p_0$ to $p_1$. 

{\rm (c)} The stabilizer of a $p$-simplex is isomorphic to $\M_{\infty,r-1+k}$ for some  $0\le k\le p$. 
\end{prop}

\begin{proof}
Let $\s,\s'$ be two $p$-simplices of $\G_\infty$. They both lie in $\G(S_{n,r};\pa_0,\pa_1)$ for some $n$ large. As in Prop.~\ref{G-action}, the 
orbits of the action of $\M_{n,r}$ on $\G(S_{n,r};\pa_0,\pa_1)$ correspond to the possible orderings of the arcs, and $\al$ induces an isomorphism 
between these sets of orbits once $n$ is large enough.

The stabilizer of $\s$ for the action of $\M_{n,r}$ on
$\G(S_{n,r};\pa_0,\pa_1)$ is isomorphic to $\M_{n-k-p,r-1+k}$ for some $0\le k\le p$. So $\M_{\infty,r-1+k}$ is a subgroup of the stabilizer of $\s$.  
On the other hand, if $g\in \M_{\infty,r}$ stabilizes $\s$, then $g$ is in $\M_{m,r}$ for some $m$ large, stabilizing $\al^{m-n}(\s)$. 
Hence $g\in\M_{m-k-p,r-1+k}$ with the same $k$ as above. 
\end{proof}

\begin{prop}\label{contrinf}
$H_i(\G_\infty/\M_\infty)=0$ for all $i>0$.
\end{prop}

The proof is identical to the proof of Lemma 3.3 of  \cite{Har85}. We give it here for completeness.

\begin{proof}
$\G_\infty/\M_\infty$ is a CW-complex with a $p$-cell for each ordering $(i_0,\dots,i_p)$ of the set $\{0,\dots,p\}$ and boundary map given by $$d_p(i_0,\dots,i_p)=\sum_{j=0}^n(-1)^j(\tau_j(i_0),\dots,\hat i_j,\dots,\tau_j(i_p))$$ where $\tau_j(n)=n$ for $n<i_j$ and $\tau_j(n)=n-1$ for $n>i_j$. Now the map $D_p:C_p(\G_\infty/\M_\infty)\to C_{p+1}(\G_\infty/\M_\infty)$ defined by 
$$D_p(i_0,\dots,i_p)=(p+1,i_0,\dots,i_p)$$ satisfies $Dd+dD={\rm Id}$, so that the identity and the zero maps are chain homotopic. 
\end{proof}

\begin{lem}\label{inject}
The map $\beta_i: H_i(\M_{n,r};\Z)\rar H_i(\M_{n,r+1};\Z)$ is injective for any $r\ge 1$ and any $n\ge 0$. 
\end{lem}

\begin{proof}
Note that $\beta:\M_{n,r}\to \M_{n+1,r}$ has a left inverse induced by gluing a disc on one of the new boundary components. 
\end{proof}

\begin{thm}\label{inf-stab}
The map $\beta_i:H_i(\M_{\infty,r};\Z)\to H_i(\M_{\infty,r+1};\Z)$ is an isomorphism for all $i$. 
\end{thm}

\begin{cor}\label{bdry}
The map $\beta_i:H_i(\M_{n,r};\Z)\rar H_i(\M_{n,r+1};\Z)$
is an isomorphism when $n\ge 4i+3$ for any $r\ge 1$.  
\end{cor}

\begin{proof}[Proof of the corollary]
This follows from the commutativity of the following diagram:
$$\xymatrix{H_i(\M_{n,r};\Z)\ar[d]_\al \ar[r]^-{\beta_i} & H_i(\M_{n,r+1};\Z) \ar[d]^\al\\
H_i(\M_{\infty,r};\Z)\ar[r]^-{\beta_i} & H_i(\M_{\infty,r+1};\Z)
}$$
where the vertical arrows are isomorphisms by  Theorem~\ref{genusstab} and the bottom one by Theorem~\ref{inf-stab}.
\end{proof}

\begin{proof}[Proof of the theorem]
The result is true when $i=0$. By the lemma, we only need to prove surjectivity. 
We do this by induction on $i$, supposing the result is true for all $i<q$.
 
We use the spectral sequences for the action of $\M_{\infty,r+1}$ on $\G_{\infty,r+1}$, i.e. for the double complex 
$$E_{q}\M_{\infty,r+1} \otimes_{\M_{\infty,r+1}} \tilde C_p(\G_{\infty,r+1}).$$
The first spectral sequence has $E^1_{p,q}=0$ as $\G_\infty$ is contractible (Thm.~\ref{Ginf}), so the second spectral sequence converges to 0. 
Its $E^1$ term is $E^1_{p,q}=\bigoplus_{\s\in\mathcal{O}_p}H_{q}(St(\s))$.
The first differential $$d^1:E^1_{0,q}=H_q(St(\s_0))\to E^1_{-1,q}=H_q(\M_{\infty,r+1})$$ is 
the map $\beta_q$ we are interested in. 
The $j$th row in the $E^1$-term of the spectral sequence computes the homology of $\G_\infty/\M_\infty$ with local coefficients in 
$H_j(St(\s))$. If $\s_p$ is a $p$-simplex, $St(\s_p)\cong \M_{\infty,r+k}$ for some $0\le k\le p$ (Prop.~\ref{action2}). 
The restriction of the differential $d^1:E^1_{p,j}\to E^1_{p-1,j}$ to a summand $H_j(St(\s))$ is given by the alternating sum of the $p+1$ 
homomorphisms corresponding to the face maps. Each face of $\s_p$ is $h\s_{p-1}$ for $\s_{p-1}$ a representative of an orbit of $(p-1)$-simplices and 
$h$ some element of $\M_{\infty,r+1}$, and the summand of $d^1$ corresponding to this face is the map $c_h:H_j(St(\s_p))\to H_j(St(\s_{p-1})$
induced by conjugation by $h$. 
These maps fit into a diagram 
$$\xymatrix{H_j(St(\s_p)) \ar[r]^-{c_h}\ar[d] & H_j(St(\s_{p-1})) \ar[d]\\
H_j(\M_{\infty,r+1}) \ar[r]^-{c_h} & H_j(\M_{\infty,r+1})
}$$ 
where the bottom map is the identity as it is given by an inner automorphism. 
By induction, when $j<q$ the vertical maps are isomorphisms and hence the coefficient system is trivial. Using Prop.~\ref{contrinf}, we get that the $E^2$-term 
is trivial below the $q$th line. So $d^1:E^1_{0,q}\to E^1_{-1,q}$ is the only map that can kill $E^1_{-1,q}$ and  
hence it must be surjective. 
\end{proof}

\section{closed surfaces}

In this section, we prove the stability theorem for closed surfaces. In a non-orientable surface, there are 4 types of embedded circles: 
there are 1-sided and 2-sided circles, with orientable or non-orientable complement. (However, 2-sided circles with orientable complement 
exist only in surfaces of even genus, and 1-sided circles with orientable complement exist only in odd genus surfaces.) 
We show that the complex of embedded circles  with connected, non-orientable complement is highly connected using the connectivity of 
the complex of all circles given in \cite{Iva87}. Part (3) of Theorem~\ref{main} is then proved
by comparing the spectral sequence for the action of the 
mapping class group on that complex for a closed surface and for a surface with one boundary component. 

\vs

We say that an embedded circle in $S$ is {\em non-degenerate} if it does not bound a disc or a M\"{o}bius band, and is not isotopic to a boundary component of $S$. 

\begin{Def} For any surface $S$, define\\
\begin{tabular}{lp{4.1in}}
$\C(S) =$& the simplicial complex whose vertices are isotopy classes of non-degenerate embedded circles in $S$. 
A $p$-simplex in $\C(S)$ is a collection of $p+1$ circles $\lgl C_0,\dots,C_p\rgl$ embeddable disjointly and non-pairwise isotopic.
\end{tabular}
\end{Def}

One could alternatively work with the complex Ivanov denotes $\bar\C(S)$ which allows the curves bounding a M\"{o}bius band. A retraction 
$\bar\C(S)\to \C(S)$ is obtained by mapping such a curve to the core of the M\"{o}bius band. 

\begin{thm}\cite[Thm. 2.6]{Iva87}\label{CS}
$\C(S)$ is $(e(S)-1)$-connected (resp. $(e(S)-2)$- or
$(e(S)-3)$-connected) when $S$ has 0 (resp. 1 or more than 2) boundary components, where $e(S)=-\chi(S)$ is the opposite of the Euler characteristic of $S$. 
\end{thm}

As for arcs, we need to consider the subcomplexes of curve systems with connected non-orientable complement. 

\begin{Def}
Let $S$ be any surface.\\
\begin{tabular}{lp{4in}}
$\C_0(S)=$& the subcomplex of $\C(S)$ consisting of simplices with  connected complement.
\end{tabular}

For $S$ non-orientable, we consider moreover the complex\\
\begin{tabular}{lp{4in}}
$\D(S) =$& the subcomplex of $\C_0(S)$ consisting of simplices with  non-orientable (connected) complement.
\end{tabular}
\end{Def}

A very natural complex to consider for a non-orientable surface $S$ is the  subcomplex of $\D(S)$ of 1-sided circles. (While 1-sided curve 
systems are always non-separating, they can have an orientable complement.) This complex is connected if the genus of $S$ is large enough, but we do not know if it is highly connected.

\begin{thm}\label{C0}  
$\C_0(S)$ is $(\frac{h-3}{2})$-connected where $h$ is the genus of $S$ or twice its genus if $S$ is orientable. 
\end{thm}

This is \cite[Thm.1.1]{Har85} when $S$ is orientable and the proof in the orientable and non-orientable case is similar to the proof 
of Theorem~\ref{BX}, the difference in connectivity bound coming from the difference between the Euler characteristic of an arc and 
that of a circle. We give a sketch of the proof for convenience.

\begin{proof}
We do an induction on the pair $(h,r)$, where $r$ is the number of boundary components of $S$. Note that the result is true for $h\le 2$ for any $r$. Fix $(h,r)\ge (3,0)$ and $k\le \frac{h-3}{2}$. 
As $\C(S)$ is at least $(h-3)$-connected (Theorem~\ref{CS}), any map $f:S^k\to \C_0(S)$ can be extended to a map $\hat f:D^{k+1}\to \C(S)$. For $\s$ a regular bad simplex of maximal dimension in $D^{k+1}$, defined as in the proof of Theorem~\ref{BX}, the restriction of $\hat f$ to the link of $\s$ maps to the join $J_\s=C_0(X_1)*\dots * \C_0(X_c)$ if $S\minus\s=X_1\sqcup\dots \sqcup X_c$. Each $\C_0(X_i)$ is $(\frac{h_i-3}{2})$-connected by induction. The Euler characteristic gives $\sum_i(2-h_i-r_i)=2-h-r$, and $\sum_i r_i=r+k$ for some $p+1\le k\le 2p+2$. The connectivity of $J_\s$ is $(\sum_i\lfloor \frac{h_i+1}{2}\rfloor)-2\ge (\frac{\sum_ih_i}{2})-2\ge \frac{h-(2p+2)+2c-2}{2}-2\ge \frac{h-4}{2}-p$ as $c\ge 2$. 
We need the connectivity of $J_\s$ to be at least $(\frac{h-3}{2}-p)$ to be able to extend $\hat f\big|_{\link(\s)}$ to a disc $D^{k+1-p}$.
If $h$ is even, 
$\lfloor\frac{h-4}{2}\rfloor=\lfloor\frac{h-3}{2}\rfloor$.  If $k<2p+2$, we gain $\frac{1}{2}$ in the next to last inequality which also gives enough connectivity. The case left is when $h$ is odd and $k=2p+2$, but in that case one of the $h_i$'s must also be odd and we gain $\frac{1}{2}$ in the first inequality above. 
The end of the proof is as in Theorem~\ref{BX}.
\end{proof}

\begin{thm}\label{closed-conn}
For $S_{n,r}$ non-orientable, $\D(S_{n,r})$ is $(\frac{n-5}{2})$-connected. 
\end{thm}

The proof is now analogous to that of Theorem~\ref{GSD} and we give a sketch for convenience.

\begin{proof}
 A map $f:S^k\to \D(S)$ can be extended to a map $\hat f:D^{k+1}\to \C_0(S)$ when $k\le \frac{n-5}{2}$ by the previous theorem. For a regular bad simplex $\s$, now defined as in the proof of Theorem~\ref{GSD}, $S\minus\hat f(\s)$ is an orientable surface of genus $g\ge \frac{n-2p-2}{2}$, and $\hat f:\link(\s)\to \C_0(S\minus\hat f(\s))$ if $\s$ is of maximal dimension. The latter space is $(\frac{n-2p-5}{2})$-connected by Theorem~\ref{C0}, so that the map can be extended over a disc $D^{k+1-p}$. Using $\hat f * F$ to modify $\hat f$ in the interior of the $\operatorname{Star}(\s)\cong \del\s * D^{k+1-p}$ improves the situation. Indeed, let $\al * \beta$ be a regular bad simplex of  $\del\s * D^{k+1-p}$. As in the proof of Theorem~\ref{GSD}, there are two cases. Either $\hat f(\al)=\hat f(\s)$ and $\al*\beta=\al$ is regular bad of smaller dimension than $\s$, or there exists a circle $C_j$ in $\hat f(\s)\minus\hat f(\al)$. As $\beta$ does not separate $S\minus\hat f(\s)$, there exists an arc in $S\minus(\hat f(\s)*F(\beta))$ joining the two copies of any point of $C_j$. Such an arc closes to a 1-sided circle in the complement of $\hat f(\al) * F(\beta)$ in $S$, so that $\hat f(\al) * F(\beta)$ has non-orientable complement, contradicting the badness assumption on $\al * \beta$. The result follows by induction.
\end{proof}

\subsection{Spectral sequence argument}

Consider the (non-augmented) spectral sequence for the action of $\M=\M_{n,r}$ on $X=\D(S_{n,r})$ with
$$E^1_{pq}=\oplus_{\s_p\in \mathcal{O}_p} H_q(St(\s_p);\Z_{\s_p})  \Rightarrow H_{p+q}^\M(X)$$
where $\mathcal{O}_p$ is a set of representatives for the orbits of $p$-simplices $\s_p$ in $X$ and $\Z_{\s_p}$ is the module $\Z$ twisted by the orientation of $\s_p$.
(The twisting in the coefficients comes from the fact that the stabilizer of a collection of circles does not need to fix the simplex pointwise.)
Following ideas of \cite{Iva93}, we prove the theorem by comparing the spectral sequence for a closed surface $S$ to that of a surface 
with one boundary component. 

Let $R$ be a surface of genus $n$ with one boundary component, and let $S$ be the surface obtained from $R$ by 
gluing a disc on its boundary. If $^RE^1_{pq}$ and $^SE^1_{pq}$ denote the two spectral sequences with $S_{n,r}$ being $R$ and $S$ respectively,  
gluing a disc induces a map of spectral sequences 
$$^RE^1_{pq} \ \rar\  ^SE^1_{pq}.$$ 

For each $p$, we have $\min(p+2,n-p-1)$ orbits of $p$-simplices, parametrized by the number of 1-sided and 2-sided circles in the simplex, and gluing 
a discs gives a 1--1 correspondence between the orbits of $p$-simplices in $\D(R)$ and those in $\D(S)$. 

Let $\s_p$ be a $p$-simplex of $\D(R)$ and $\s_p'$ the corresponding simplex in $\D(S)$. The stabilizer $St_R(\s_p)$ is a subgroup of 
$\M^+_{n-p-k-1,p+k+2}$ where $0\le k\le p+1$
is the number of $2$-sided circles in $\s_p$ and where $\M^+$ denotes the extended mapping class group where the boundaries of the surface 
are not fixed pointwise and may be permuted. 
Similarly, $St_S(\s_p')\le \M^+_{n-p-k-1,p+k+1}$. Gluing a disc on the boundary of $S$ induces a map 
$H_i(\M_{n-p-k-1,p+k+2})\to H_i(\M_{n-p-k-1,p+k+1})$ which is always surjective (because it is a right inverse to gluing a pair of pants) and is an isomorphism 
when $n-p-k-1\ge 4i+3$ by Corollary~\ref{bdry}.
We need such an isomorphism on the homology of the stabilizers to obtain an isomorphism of spectral sequences.

\begin{lem}
The map $H_i(St_R(\s_p);\Z_{\s_p})\rar H_i(St_S(\s_p');\Z_{\s_p'})$ is an isomorphism when $i\le \frac{n-2p-5}{4}$. 
Moreover, the map is surjective when $i\le \frac{n-2p-1}{4}$.
\end{lem}

\begin{thm}\label{closedstab}
The map $H_i(\M_{n,1};\Z)\rar H_i(\M_{n,0};\Z)$ is an isomorphism when $n\ge 4i+5$. Moreover, the map is surjective when $n\ge 4i+1$. 
\end{thm}

\begin{proof}[Proof of the Theorem]
We only need to prove the result when $i\ge 1$.
By the Lemma, we have an isomorphism of spectral sequences in the range $q\le \frac{n-2p-5}{4}$, so in particular when $p+q\le \frac{n-5}{4}$, and 
an epimorphism when $p+q\le \frac{n-1}{4}$. Then, by \cite[Thm.~1.2]{Iva93}, we have an isomorphism 
$H_i^\M(\D(R))\to H_i^\M(\D(S))$ when $i\le \frac{n-5}{4}$ and an epimorphism when $i\le \frac{n-1}{4}$. 
As $i\le \frac{n-5}{4}$ (with $i\ge 1$) implies $i\le \frac{n-5}{2}$, the isomorphism part of the result follows from Theorem~\ref{closed-conn}. For the surjectivity, we need moreover that $i\le \frac{n-1}{4}$ implies   $i\le \frac{n-5}{2}$. This is only true when $n\ge 9$. The surjectivity however holds without this extra assumption since for $n\le 8$ it only concernes $H_1$, and surjectivity always holds on $H_1$ because the homomorphism $\M_{n,1}\to \M_{n,0}$ is surjective.
\end{proof}

\begin{proof}[Proof of the Lemma]
Let $\widetilde{St}_R(\s_p)\le St_R(\s_p)$ and $\widetilde{St}_S(\s_p')\le St_S(\s'_p)$ denote the subgroups of elements 
fixing each circle of $\s_p$ or $\s_p'$, preserving its orientation. 
Let $\Si_k^\pm$ denote the group of signed permutations of $k$ elements. Suppose that $\s_p$ has exactly $k$ 2-sided circles. We have a 
commutative diagram with exact rows:
$$\xymatrix{1\ar[r] & \widetilde{St}_R(\s_p) \ar[r] \ar[d] & St_R(\s_p) \ar[r] \ar[d] & \Si_k^\pm \times \Si_{p+1-k}^\pm \ar[r] \ar[d]^= & 1\\
1\ar[r] & \widetilde{St}_S(\s_p') \ar[r] & St_S(\s_p') \ar[r] & \Si_k^\pm \times \Si_{p+1-k}^\pm \ar[r]  & 1
}$$
Each row gives rise to a Hochschild-Serre spectral sequence with 
$$^RE^2_{rs}=H_r(\Si_k^\pm \times \Si_{p+1-k}^\pm;H_s( \widetilde{St}_R(\s_p);\Z_{\s_p})) \Rightarrow H_{r+s}(St_R(\s_p);\Z_{\s_p})$$
and $E^1$-term $E^1_{rs}=F_r\otimes_{\Si_k^\pm\times\Si_{p+1-k}^\pm}H_s( \widetilde{St}_R(\s_p);\Z_{\s_p})$ where $F_r$ is a projective resolution 
of $\Z$ over $\Z[\Si_k^\pm\times\Si_{p+1-k}^\pm]$,  
 and similarly for the bottom row. We have a map of spectral sequences $\psi:\, ^RE^1_{rs}\to\,^SE^1_{rs}$. If we can show that 
$H_s( \widetilde{St}_R(\s_p);\Z_{\s_p})\to H_s( \widetilde{St}_R(\s_p);\Z_{\s_p})$ is an isomorphism whenever $s\le \frac{n-2p-5}{4}$, the map $\psi$ 
will be an isomorphism of 
spectral sequences in this range for any $r$, so in particular when $r+s\le \frac{n-2p-5}{4}$ and when $r+s\le \frac{n-2p-1}{4}$ if $r\ge 1$.  
The first part of the lemma will then follow by \cite[Thm.~1.2]{Iva93}. For the surjectivity part, we need moreover the map to be 
surjective when $s\le \frac{n-2p-1}{4}$.

Note that the action of  $\widetilde{St}_R(\s_p)$ on the twisted coefficient $\Z_{\s_p}$ is actually trivial, so it can be 
replaced by $\Z$. 

There is a second commutative diagram with exact rows:
$$\xymatrix{1\ar[r] & \Z^{p+1} \ar[r] \ar[d] & \M(R\minus\s_p) \ar[r] \ar[d] & \widetilde{St}_R(\s_p) \ar[r] \ar[d] & 1\\
1\ar[r] & \Z^{p+1} \ar[r] & \M(S\minus\s_p') \ar[r] & \widetilde{St}_S(\s_p') \ar[r] & 1
}$$
where the first $k$ factors in $\Z^{p+1}$ are generated by Dehn twists of the form $t_{C_i^+}t^{-1}_{C_i^-}$ with $C_i^+$ and $C_i^-$ the two 
boundary components in $S\minus\s_p$ coming from a 2-sided circle $C_i$ of $\s$, and the other factors come from Dehn twists along boundary components coming 
from 1-sided circles. Also, $\M(R\minus\s_p)$ is isomorphic to $\M_{n-p-k-1, p+k+2}$ and  $\M(S\minus\s_p')$ is isomorphic to $\M_{n-p-k-1, p+k+1}$.
We again have a map of Hochschild-Serre spectral sequences  $^RE^2_{rs}\to\, ^SE^2_{rs}$ from the diagram, where now 
$$^RE_{r,s}^2=H_r(\widetilde{St}_R(\s_p);H_s(\Z^{p+1})) \Rightarrow H_{r+s}(\M(R\minus\s_p)).$$ 
We have $^RE^2_{0,s}\sta{\cong}{\to}\, ^SE^2_{0,s}$ for all $s$ because $\Z^{p+1}$ is central in 
both groups. We know moreover that the sequences converge to isomorphic groups when $n-2p-2\ge 4(r+s)+3$ by 
Corollary~\ref{bdry} as $k\le p+1$.  
(The map, induced by gluing a disc, is the right inverse to gluing a pair of pants in this case.)
By \cite[Thm.~1.3]{Iva93} (using the universal coefficient theorem for condition (v)), we get an isomorphism $^RE^2_{r,0}\to\, ^SE^2_{r,0}$
for all $r\le \frac{n-2p-5}{4}$ which is the condition we needed above for the isomorphism in the lemma. For the surjectivity when $r\le \frac{n-2p-1}{4}$, we need moreover 
that  $H_{r+s}(\M(R\minus\s_p))\to H_{r+s}(\M(S\minus\s_p))$ is surjective when $r+s \le \frac{n-2p-1}{4}$ and this map is always surjective. 
\end{proof}

\section{Stable homology}\label{homo}

In this section, we describe the stable homology of the non-orientable mapping class groups. 
We first relate the non-oriented (1+1)-cobordism category to the stable non-orientable mapping class group, 
using a combinatorial model of this category. 
We then use the work of Galatius-Madsen-Tillmann-Weiss which computes the homotopy type of the cobordism category in 
terms of a Thom spectrum associated to the Grassmannians of 2-planes in $\RR^{n+2}$.

\vs

\begin{figure}
\includegraphics{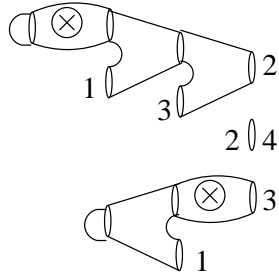}   
\caption{Example of 1-morphism in $\CCC_{2,b}$ from 3 to 4}\label{catC}
\end{figure}

Let $\CCC_2$ denote the $(1+1)$-cobordism 2-category, namely $\CCC_2$ has objects the natural numbers, where $n$ is thought of as a disjoint 
union of $n$ circles, 1-morphisms from $m$ to $n$ are cobordisms with $m$ incoming and $n$ outgoing circles, and 2-morphisms are given by 
isotopy classes of diffeomorphisms of the cobordisms fixing the boundary pointwise. 
We denote by $\CCC_{2,b}$ the subcategory where the natural map $\{1,\dots,n\}\to \pi_0(S)$ is surjective for any 1-morphism $S$ with $n$ 
outgoing circles. 
 More precisely, we use the following combinatorial model for $\CCC_{2,b}$, which is a variation of Tillmann's model for the oriented cobordism category in \cite{Til99}: fix a pair of pants $P$, a punctured M\"{o}bius band $M$ 
 and a disc $D$ which are respectively 1-morphisms 
$$P:2 \to 1,\ \ \ \bar P:1\to 2,\ \ \   M:1\to 1\ \ \ \textrm{ and }\ \ \ D:0\to 1.$$  
The boundary circles of these surfaces are parametrized by $[0,2\pi[$. 
We also consider the circle $C$ as a 1-morphism from 1 to 1. 
A 1-morphism in $\CCC_{2,b}$ from $m$ to $n$ is then a surface build out of $P,\bar P,M,D$ and $C$ by gluing and disjoint union, with the incoming 
boundary components labeled $1,\dots,m$ and the outgoing ones labeled $1,\dots,n$. (See Fig.~\ref{catC} for an example.)
Isolated circles are labeled on both sides. 
Note that such a morphism has at most $n$ components. In particular, there is no morphism to 0. 
The identity 
morphism from $n$ to $n$ is a disjoint union of $n$ circles with the same labels on both sides. The other labellings of the circles give an inclusion of the 
symmetric group $\Si_n$ in $\CCC_{2,b}(n,n)$ which acts on other morphisms by permuting the labels. Composition of 1-morphisms is by gluing.

The 2-morphisms are isotopy classes of diffeomorphisms taking one surface to the other, identifying the boundary circles via the identity 
map according to the labels. 

Following the notation of \cite{Til97}, let $\BB\CCC_{2,b}$ denote the category enriched over simplicial sets obtained 
from $\CCC_{2,b}$ by taking the classifying spaces (or nerve) of the categories of morphisms,  
i.e. $\BB\CCC_{2,b}(m,n):=\B\,(\CCC_{2,b}(m,n))$. 
As 1-morphisms from $k$ to $1$ in $\CCC_{2,b}$ are connected surfaces with $k+1$ boundary components, we have 
$$\BB\CCC_{2,b}(k,1)\simeq ({\textstyle\coprod_{g\ge 0}}\,\B\Ga_{g,k+1})\ {\textstyle \coprod}\ ({\textstyle\coprod}_{n\ge 0}\,\B\M_{n,k+1})$$
where $\Ga_{g,k+1}$ denotes the mapping class group of an orientable surface of genus $g$ with $k+1$ boundary components with  the boundaries  
fixed pointwise. (As we are only considering surfaces with at least one boundary component here, the diffeomorphisms automatically preserve the orientation.) 
We denote by $\B\CCC_{2,b}$ the classifying space of the simplicial category $\BB\CCC_{2,b}$. 

It is equivalent to work with the space of all diffeomorphisms instead of just isotopy classes 
as the components of the diffeomorphism groups are contractible 
except in a few low genus cases
\cite{EarEel,EarSch} and we will let the genus tend to infinity. 
In particular, the category $\BB\CCC_2$ can be replaced in the following theorem 
by the category of embedded cobordisms described in the introduction.

\begin{thm}\label{til}
$H_*(\Om\B\CCC_{2,b};\Z) \cong H_*(\Z\x \B\M_\infty;\Z)$. 
\end{thm}

This is the {\em  non-oriented} analogue of \cite[Thm.~3.1]{Til97}.  
Note indeed that $\CCC_{2,b}$ is the cobordism category of all surfaces, not just the non-orientable ones. 
We give a sketch of the proof of the theorem for convenience. It is completely analogous to the proof of \cite[Thm.~3.1]{Til97}
but relies on our stability theorem. 

\begin{proof}
Let $\CCC_\infty$ be the $\BB\CCC_{2,b}$-diagram defined by 
$$\CCC_\infty(k)=\colim(\BB\CCC_{2,b}(k,1)\sta{\circ M}{\rar} \BB\CCC_{2,b}(k,1)\sta{\circ M}{\rar}\dots).$$ 
The maps $\BB\CCC_{2,b}(m,k)\times \CCC_\infty(k)\to \CCC_\infty(m)$ are induced by pre-composition in $\BB\CCC_{2,b}$.  
We have $\CCC_\infty(k)\simeq \Z\times \B\M_{\infty,k+1}$.
Let $E_{\CCC_{2,b}}\CCC_\infty={\rm hocolim}_{\BB\CCC_{2,b}}\CCC_\infty$. 
This is a contractible space, being a colimit of contractible spaces (see \cite[Lem.~3.3]{Til97}).
For every object $k$, there is a pull-back diagram 
$$\xymatrix{\CCC_\infty(k) \ar[r] \ar[d] & E_{\BB\CCC_{2,b}}\CCC_\infty \ar[d]\\
k \ar[r] & \B\CCC_{2,b}
}$$
For each vertex in $v\in \B\CCC_{2,b}(m,k)$, the induced map $v_*:H_*(\CCC_\infty(k),\Z)\to H_*(\CCC_\infty(m),\Z)$ is an isomorphism by Theorem~\ref{main}. 
Indeed, it is enough to show that each building block $P,\bar P,M$ and $D$ induces an isomorphism in homology. Theorem~\ref{main} implies that $P,M$ and $D$ induce
isomorphims. It follows that $\bar P$ also induces an isomorphism as $P\circ\bar P\circ M\cong M\circ M\circ M$.   

It follows from \cite[Thm.~3.2]{Til97} that the above diagram is homology Cartesian, i.e. that the inclusion of the fiber into the homotopy fiber 
is a homology isomorphism. And the homotopy fiber is $\Om \B\CCC_{2,b}$ as the total space is contractible.
\end{proof}

\vs

Let $G_n=Gr(2,n+2)$ be the Grassmannian of (non-oriented) 2-planes in $\RR^{n+2}$. Let $U_n$ and $U_n^\perp$ denote the tautological 
bundle over $G_n$ and its orthogonal complement, and let $Th(U_n^\perp)$ denote the Thom space of $U_n^\perp$, the disc bundle of $U_n^\perp$ 
with its sphere bundle collapsed to a point. 
Pulling back $U^\perp_{n+1}$ along the inclusion $G_n\inc G_{n+1}$ gives a map 
$$S^1\wedge Th(U_n^\perp)\to Th(U_{n+1}^\perp)$$
so that the spaces $Th(U_n^\perp)$ form a spectrum. Following \cite{GMTW05}, we denote $\GG_{-2}$ the spectrum with $(\GG_{-2})_n=Th(U^\perp_{n-2})$. 
Let $$\Om^\infty \mathbb{G}_{-2}:=\colim_{n\to\infty} \Om^{n+2}Th(U_n^\perp)$$
be its associated infinite loop space. 
The main theorem of \cite{GMTW05}, in the (non-oriented) dimension 2 case, says  that 
$$ \B\CCC_2 \simeq \Om^\infty \Si \mathbb{G}_{-2}$$
and thus that $\Om\B\CCC_2\simeq \Om^\infty\GG_{-2}$.  (They use the embedded cobordisms model for $\CCC_2$.)
This homotopy equivalence is induced by the 
Thom-Pontrjagin construction explained in the introduction. Furthermore, they show that 
$$\B\CCC_2\simeq \B\CCC_{2,b}.$$

Combined with our Theorem~\ref{til}, this yields

\begin{thm}\label{integral}
 $H_*(\M_\infty;\Z)\cong H_*(\Om_0^\infty \mathbb{G}_{-2};\Z)$
\end{thm}

Here, $\Om_0^\infty \mathbb{G}_{-2}$ denotes the 0th component of $\Om^\infty \mathbb{G}_{-2}$.

Theorem \ref{integral} can alternatively be proved running the original proof of the Madsen-Weiss Theorem \cite{MadWei02} removing orientations everywhere 
and replacing Harer's stability theorem by our Theorem~\ref{main}.

\begin{cor}\label{awayfrom2}
 $H_*(\M_\infty;\Z[\frac{1}{2}])\cong H_*(\Om^\infty_0\Si^\infty (\BO(2)_+);\Z[\frac{1}{2}])$
\end{cor}

\begin{cor}\label{rational}
 $H^*(\M_\infty;\QQ)\cong \QQ[\zeta_1,\zeta_2,\dots]$ with $|\zeta_i|=4i$.
\end{cor}

\begin{proof}[Proof of Corollary~\ref{awayfrom2} and \ref{rational}]
The cofiber sequence of spectra 
$$Th(U_n^\perp\big|_{S(U_n)}) \rar Th(U_n^\perp) \rar Th(U_n\oplus U_n^\perp)$$
yields the following homotopy fibration sequence of infinite loop spaces 
$$\Om^\infty \GG_{-2} \rar \Om^\infty\Si^\infty(\BO(2)_+) \rar \Om^\infty\RPi_{-1} $$
where $\Om^\infty\RPi_{-1}=\Om^\infty\GG_{-1}$ is the infinite loop space associated to the Thom spectrum obtained as above from the orthogonal 
bundle to the canonical bundle over the Grassmannian of lines in $\RR^{n+1}$ (see \cite[Prop.~3.1]{GMTW05}).

The analogous cofiber sequence of spectra for $\RPi_{-1}$ yields the fibration sequence:
$$\Om^\infty\RPi_{-1}  \rar \Om^\infty\Si^\infty(\RPi_+) \sta{\del}{\rar} \Om^\infty S^\infty$$
where $\del$ is the stable transfer associated with the universal double covering space $p:\operatorname{E}\Z/2 \to \B\Z/2$ \cite[Rem.~3.2]{GMTW05}. 
As $p$ is an equivalence away from 2, so is $\del$ by \cite[Thm.~5.5]{BecGot}, and the homology of $\Om^\infty\RPi_{-1}$ must be 2-torsion, which gives Corollary~\ref{awayfrom2}. 
Finally for Corollary~\ref{rational} 
note first that $H^*(\Om^\infty_0\Si^\infty(\BO(2)_+),\QQ)\cong H^*(\Om^\infty_0\Si^\infty(\BO(2)),\QQ)$ as $\Om^\infty\Si^\infty(\BO(2)_+)\simeq\Om^\infty\Si^\infty(\BO(2))\times\Om^\infty S^\infty$ and $\Om^\infty S^\infty$ has trivial rational cohomology. 
From the theory of graded commutative Hopf algebras \cite[Append.]{MilMoo65}, we have that $H^*(\Om^\infty_0\Si^\infty(\BO(2)),\QQ)$ is the symmetric algebra on $H^{>0}(\BO(2),\QQ)$. (See also \cite[Thm.~2.10]{Mad06}, and \cite[Thm.~3.1]{GMT06} for identifying the classes as the classes described in the introduction.)
\end{proof}

Recall from \cite[Thm.~5.1]{DyeLas62}\cite{CohLadMay} that for any space $X$, the homology of $\Om^\infty\Si^\infty X$ at any prime $p$ 
is an algebra over the homology of $X$, 
generated by the Dyer-Lashof operations. We know the homology of $\BO(2)$  \cite{Bro82,Fes83} and of $\RPi$ so that, 
using the two fibration sequences above, it is possible to compute the homology of $\Om^\infty\GG_{-2}$, and hence that 
of $\M_\infty$. The analogous computations in the orientable and spin cases were carried out  
by Galatius in \cite{Gal04,Galspin}.
One can verify as a first step that $H_1(\Om^\infty\GG_{-2};\Z)=\pi_1(\GG_{-2})=\Z/2$.

\bibliography{biblio}

\begin{thebibliography}{10}

\bibitem{BecGot}
J.~C. Becker and D.~H. Gottlieb.
\newblock The transfer map and fiber bundles.
\newblock {\em Topology}, 14:1--12, 1975.

\bibitem{Bro82}
Edgar~H. Brown, Jr.
\newblock The cohomology of {$B{\rm SO}\sb{n}$} and {$B{\rm O}\sb{n}$} with
  integer coefficients.
\newblock {\em Proc. Amer. Math. Soc.}, 85(2):283--288, 1982.

\bibitem{CohLadMay}
Frederick~R. Cohen, Thomas~J. Lada, and J.~Peter May.
\newblock {\em The homology of iterated loop spaces}.
\newblock Springer-Verlag, Berlin, 1976.
\newblock Lecture Notes in Mathematics, Vol. 533.

\bibitem{DyeLas62}
Eldon Dyer and R.~K. Lashof.
\newblock Homology of iterated loop spaces.
\newblock {\em Amer. J. Math.}, 84:35--88, 1962.

\bibitem{EarEel}
C.~J. Earle and James Eells.
\newblock A fibre bundle description of {T}eichm\"uller theory.
\newblock {\em J. Differential Geometry}, 3:19--43, 1969.

\bibitem{EarSch}
C.~J. Earle and A.~Schatz.
\newblock Teichm\"uller theory for surfaces with boundary.
\newblock {\em J. Differential Geometry}, 4:169--185, 1970.

\bibitem{Fes83}
Mark Feshbach.
\newblock The integral cohomology rings of the classifying spaces of {${\rm
  O}(n)$} and {${\rm SO}(n)$}.
\newblock {\em Indiana Univ. Math. J.}, 32(4):511--516, 1983.

\bibitem{Gal04}
S{\o}ren Galatius.
\newblock Mod {$p$} homology of the stable mapping class group.
\newblock {\em Topology}, 43(5):1105--1132, 2004.

\bibitem{Galspin}
S{\o}ren Galatius.
\newblock Mod 2 homology of the stable spin mapping class group.
\newblock {\em Math. Ann.}, 334(2):439--455, 2006.

\bibitem{GMT06}
Soren Galatius, Ib~Madsen, and Ulrike Tillmann.
\newblock Divisibility of the stable {M}iller-{M}orita-{M}umford classes.
\newblock {\em J. Amer. Math. Soc.}, 19(4):759--779 (electronic), 2006.

\bibitem{GMTW05}
Soren Galatius, Ib~Madsen, Ulrike Tillmann, and Michael Weiss.
\newblock {The homotopy type of the cobordism category}.
\newblock arXiv:math.AT/0605249.

\bibitem{Gra73}
Andr{\'e} Gramain.
\newblock Le type d'homotopie du groupe des diff\'eomorphismes d'une surface
  compacte.
\newblock {\em Ann. Sci. \'Ecole Norm. Sup. (4)}, 6:53--66, 1973.

\bibitem{Han07}
Elizabeth Handbury.
\newblock {Homology stability of non-orientable mapping class groups with
  marked points}.
\newblock preprint.

\bibitem{Har85}
John~L. Harer.
\newblock Stability of the homology of the mapping class groups of orientable
  surfaces.
\newblock {\em Ann. of Math. (2)}, 121(2):215--249, 1985.

\bibitem{Hat91}
Allen Hatcher.
\newblock On triangulations of surfaces.
\newblock {\em Topology Appl.}, 40(2):189--194, 1991.

\bibitem{HatVogWah}
Allen Hatcher, Karen Vogtmann, and Nathalie Wahl.
\newblock Erratum to: ``{H}omology stability for outer automorphism groups of
  free groups" [{A}lgebr. {G}eom. {T}opol. {4}:1253--1272, 2004] by {H}atcher
  and {V}ogtmann.
\newblock {\em Algebr. Geom. Topol.}, 6:573--579, 2006.

\bibitem{HatWah07}
Allen Hatcher and Nathalie Wahl.
\newblock {Stabilization for mapping class groups of 3-manifolds}.
\newblock preprint.

\bibitem{Iva87}
N.~V. Ivanov.
\newblock Complexes of curves and {T}eichm\"uller modular groups.
\newblock {\em Uspekhi Mat. Nauk}, 42(3(255)):49--91, 255, 1987.

\bibitem{Iva90}
N.~V. Ivanov.
\newblock Stabilization of the homology of {T}eichm\"uller modular groups.
\newblock {\em Algebra i Analiz}, 1(3):110--126, 1989.

\bibitem{Iva93}
Nikolai~V. Ivanov.
\newblock On the homology stability for {T}eichm\"uller modular groups: closed
  surfaces and twisted coefficients.
\newblock In {\em Mapping class groups and moduli spaces of Riemann surfaces
  (G\"ottingen, 1991/Seattle, WA, 1991)}, volume 150 of {\em Contemp. Math.},
  pages 149--194. Amer. Math. Soc., Providence, RI, 1993.

\bibitem{Kor98}
Mustafa Korkmaz.
\newblock First homology group of mapping class groups of nonorientable
  surfaces.
\newblock {\em Math. Proc. Cambridge Philos. Soc.}, 123(3):487--499, 1998.

\bibitem{Mad06}
Ib~Madsen.
\newblock {Moduli spaces from a topological viewpoint}.
\newblock to appear in the Proceedings of the International Congress of
  Mathematicians, Madrid, Spain, 2006.

\bibitem{MadTil01}
Ib~Madsen and Ulrike Tillmann.
\newblock The stable mapping class group and {$Q(\mathbb{C} P^\infty_+)$}.
\newblock {\em Invent. Math.}, 145(3):509--544, 2001.

\bibitem{MadWei02}
Ib~Madsen and Michael~S. Weiss.
\newblock The stable moduli space of {R}iemann surfaces: Mumford's conjecture.
\newblock {\em Ann. of Math.}, 165(3):843--941, 2007.

\bibitem{Mil86}
Edward~Y. Miller.
\newblock The homology of the mapping class group.
\newblock {\em J. Differential Geom.}, 24(1):1--14, 1986.

\bibitem{MilMoo65}
John~W. Milnor and John~C. Moore.
\newblock On the structure of {H}opf algebras.
\newblock {\em Ann. of Math. (2)}, 81:211--264, 1965.

\bibitem{Mor87}
Shigeyuki Morita.
\newblock Characteristic classes of surface bundles.
\newblock {\em Invent. Math.}, 90(3):551--577, 1987.

\bibitem{Mum83}
David Mumford.
\newblock Towards an enumerative geometry of the moduli space of curves.
\newblock In {\em Arithmetic and geometry, Vol. II}, volume~36 of {\em Progr.
  Math.}, pages 271--328. Birkh\"auser Boston, Boston, MA, 1983.

\bibitem{Stu07}
Michal Stukow.
\newblock {Generating mapping class groups of non-orientable surfaces with
  boundaries}.
\newblock arXiv:math.GT/0707.3497.

\bibitem{Til97}
Ulrike Tillmann.
\newblock On the homotopy of the stable mapping class group.
\newblock {\em Invent. Math.}, 130(2):257--275, 1997.

\bibitem{Til99}
Ulrike Tillmann.
\newblock A splitting for the stable mapping class group.
\newblock {\em Math. Proc. Cambridge Philos. Soc.}, 127(1):55--65, 1999.

\bibitem{Vog79}
Karen Vogtmann.
\newblock Homology stability for {${\rm O}\sb{n,n}$}.
\newblock {\em Comm. Algebra}, 7(1):9--38, 1979.

\bibitem{Zaw04}
Myint Zaw.
\newblock The homology groups of moduli spaces of {K}lein surfaces with one
  boundary curve.
\newblock {\em Math. Proc. Cambridge Philos. Soc.}, 136(3):599--615, 2004.

\end{thebibliography}

\end{document}